\documentclass[12pt]{article}
\usepackage{amssymb}
\usepackage{amssymb,amsfonts}
\usepackage{amsmath}
\usepackage{cases}

\usepackage{indentfirst}
\usepackage[mathscr]{eucal}
\usepackage[OT1]{fontenc}

\usepackage[round,sort&compress]{natbib}
\usepackage{lscape}
\usepackage{multirow}
\usepackage{graphicx}
\usepackage{graphics}
\usepackage{color,xcolor}
\usepackage{float}

\usepackage{bm}
\usepackage{enumitem}
\usepackage{caption} 

\usepackage[colorlinks,linkcolor=blue,anchorcolor=blue,citecolor=blue]{hyperref}

\usepackage{epstopdf}
\usepackage{booktabs}

\usepackage{tabularx}
\usepackage{caption}
\usepackage{subcaption}
\usepackage{mathrsfs}
\usepackage{amsmath}
\allowdisplaybreaks[4]

\DeclareSymbolFont{largesymbol}
{OMX}{yhex}{m}{n}
\DeclareMathAccent{\Widehat}
{\mathord}{largesymbol}{"62}

\usepackage{geometry}
\newtheorem{theorem}{Theorem}

\newtheorem{remark}{Remark}[section]

\usepackage[title]{appendix}
\usepackage{appendix}

\usepackage{titlesec}

\numberwithin{equation}{section}

\def\I{{\bf {I}}}
\def\1{{\bf {1}}}
\def\0{{\bf {0}}}

\def\N{\mathcal{N}}
\def\x{{\bf {x}}}
\def\P{{\bf {P}}}

\def\W{{\bf {W}}}

\def\R{{\mathcal{R}}}

\newcommand{\bSigma}{\mbox{$\bm\Sigma$}}
\newcommand{\bmu}{\mbox{$\bm\mu$}}

\title{Linear hypothesis testing in high-dimensional heteroscedastics via random integration}
\author{Mingxiang Cao$^{*}$, Hongwei Zhang, Kai Xu \& Daojiang He\\
    \it\footnotesize  Department of Statistics, Anhui Normal University, Wuhu, China}
\date{}

\begin{document}

\maketitle \footnotetext[1]{*Corresponding author.}

\footnotetext{E-mail:\url{caomingx@163.com}}

\noindent\textbf{Abstract}\ \ In this paper, for the problem of heteroskedastic general linear hypothesis testing (GLHT) in high-dimensional settings, we propose a random integration method based on the reference L$^{2}$-norm to deal with such problems. The asymptotic properties of the test statistic can be obtained under the null hypothesis when the relationship between data dimensions and sample size is not specified. The results show that it is more advisable to approximate the null distribution of the test using the distribution of the chi-square type mixture, and it is shown through some numerical simulations and real data analysis that our proposed test is powerful.



\vspace{1cm}

\section{Introduction}\label{sec1}
\noindent 

In recent years, with the development of science and technology and social progress, data has become an indispensable part of life, and the amount of data is getting bigger and bigger. With the advent of the big data era, high-dimensional data has emerged, characterized by data with dimension p much larger than the sample size n. The analysis of high-dimensional data is an important area in statistics and machine learning, which deals with data sets with many features or variables. High-dimensional data is widely available in various fields, including gene expression, finance, etc. However, how to analyze high-dimensional data has brought great challenges to statisticians. For example, in the study of gene expression, many existing classical tests tend to be performed with their dimensions fixed and sample sizes tending to infinity under the circumstances of the test. At this time, the classical test method may become less efficient. The analysis of high-dimensional data for the statistician has brought great challenges. With the increase in sample size, statisticians began to study the high-dimensional k-sample mean problem. The k-sample mean test includes the k-sample Behrens-Fisher problem and the linear hypothesis testing problem.

For the k-sample problem, denote the sample size of k independent p-dimensional  i.i.d. samples by $n_{1},\ldots,n_{k}$. Suppose $ \mathbf{y}_{i1},\ldots,\mathbf{y}_{in_{i}}$ are i.i.d. random samples with the mean vector and covariance matrix $\bmu_{i}$ and $\bSigma_{i}$, $i\in\{1,\ldots\,k\}$, respectively. We consider the following hypothesis:
\begin{eqnarray}\label{H:1.1}
H_0: \mu_{1}=\mu_{2}=\cdots=\mu_{k}~~~~vs.~~~~~H_1: \mbox{not}\ H_0.
\end{eqnarray}
for the hypothesis in (\ref{H:1.1}), we do not need to assume equality of the k sample covariances and allow the sample dimension p to be much larger than the sample size $n=\sum_{i=1}^kn_{i}$. When k = 2, the problem (\ref{H:1.1}) reduces to the high-dimensional two-sample problem. \cite{Bai and Saranadasa:1996} derived the asymptotic properties of the classical Hotelling's T$^{2}$ test and Dempster's inexact test for two-sample problems. \cite{Schott:2007} constructed a multivariate variance test by revising \cite{Bai and Saranadasa:1996} two-sample test. \cite{Yamada and Himeno:2015} and \cite{Hu et al:2017} extended \cite{Chen et al:2010} two-sample test based on the U-statistic to the high-dimensional k-sample problem. \cite{Aoshima et al:2015} investigated the high-dimensional k-sample test problem under weaker conditions and tested the asymptotic properties of the statistic. \cite{Chen et al:2019} proposed an efficient sparse weak mean difference test, and \cite{Hiroki et al:2020} studied the two-way MANOVA problem. \cite{Zhang et al:2022(a)} studied the general linear hypothesis testing (GLHT) problem in heteroskedastic one-way MANOVA for high-dimensional data and proposed a test based on the L$^{2}$-norm of the normal distribution.

In this paper, we want to test the following heteroscedastic general linear hypothesis testing (GLHT) problem:
\begin{eqnarray}\label{H:1.2}
H_0: \tilde{G}M=0~~~~vs.~~~~~H_1:  \tilde{G}M\neq0,
\end{eqnarray}
where $\tilde{G}$ is a q$\times$k known coefficient matrix with full row rank q$<$k and $M=(\mu_{1}, \ldots,\mu_{k})^{T}$ is a k$\times$p comprised of the k mean vectors. The hypothesis in (\ref{H:1.2}) is very general. It contains some special hypothesis testing problems when the coefficient matrix $\tilde{G}$ is set differently. For example, the hypothesis in (\ref{H:1.2}) simplifies to the one-way MANOVA problem (\ref{H:1.1}) by defining $\tilde{G}$ to be any (k$-$1)$\times$k contrast matrix,i.e., which means any (k$-$1)$\times$k matrix with linearly independent rows and zero row sums. Furthermore, a variety of post-hoc and contrast tests can be expressed as (\ref{H:1.2}). To test if $-\mu_{1}+2\mu_{2}-3\mu_{3}=0$, for instance, let $e_{r,l}$ represent a unit vector of length l with the rth item being 1 and the others 0, and set $\tilde{G}=(-e_{1,k}+2e_{2,k}-3e_{3,k})^{T}$. When we set $\tilde{G}$ to be a k-dimensional row vector $(g_{11},\ldots,g_{1k})$, it results in the following hypothesis testing problem on linear combinations of k means:
\begin{eqnarray}\label{H:1.3}
H_0:\sum_{i=1}^kg_{1i}\mu_{i}=0~~~~vs.~~~~~H_1: \sum_{i=1}^kg_{1i}\mu_{i}\neq0,
\end{eqnarray}
where $i\in\{1,\ldots\,k\}$.

This assumption is a special case of the GLHT problem; \cite{Nishiyama et al:2013} examined this particular hypothesis testing problem.The GLHT problem has been studied by several scholars. In the context of multiple linear regression modeling, \cite{Fujikoshi et al:2004} considered the test and investigated Dempster's test. \cite{Zhang et al:2017} considered the GLHT problem with a common covariance matrix and proposed a test based on the L$^{2}$-norm. \cite{Zhou et al:2017} proposed a test based on the k sample U-statistic in the context of heteroscedasticity ensembles. \cite{Zhang et al:2022(b)} proposed and investigated a test statistic for the GLHT problem based on L$^{2}$-norm and constructed an adaptive test using Box $\mathrm{\chi^{2}}$-approximation.Recently, \cite{Jiang:2022} studied the two-sample mean test problem based on the random integration method and constructed a test statistic. The test is superior to existing tests in many cases and requires fewer overall parameters. Due to the advantages of the random integration method for constructing tests, this paper will study the GLHT problem based on this method.

The rest of the paper is organized as follows: In Section 2, we use random integration techniques to propose the statistic and obtain its asymptotic properties. Section 3 conducts simulations to assess the proposed test's performance on a finite sample. In Section 4, a real dataset is examined in order to contrast the proposed test with some existing methods. In Section 5, we make a few final observations. The technical proofs of the main theorems are arranged in the Appendix and some additional simulation results are provided as supplementary materials.

\section{Test procedure and main results}\label{sec2}
\subsection{Test procedure}\label{sec2.1}
\noindent 
First, the following transformations are applied to the coefficient matrix:
\begin{align}\label{eq2.1}
\tilde{G}\rightarrow\P\tilde{G},
\end{align}
after the above transformations, the GLHT problem (\ref{H:1.2}) is invariant. Where P denotes any q$\times$q non-singular matrix. Thus, our proposed test is also invariant when $\tilde{G}$ is subjected to non-singular transformations. On this basis, we can rewrite the GLHT problem (\ref{H:1.2}) as:
\begin{eqnarray}\label{H:2.1}
H_0: GM=0~~~~vs.~~~~~H_1:  GM\neq0,
\end{eqnarray}
where $G=(\tilde{G}D\tilde{G}^{T})^{-1/2}\tilde{G}, D=diag(1/n_{1},\ldots,1/n_{k}).$~Similar to \cite{Zhang et al:2022(b)}, the GLHT problem (\ref{H:2.1}) can be equivalently transformed into the following question:
\begin{eqnarray}\label{H:2.2}
H_0: C\mu=0~~~~vs.~~~~~H_1:  C\mu\neq0,
\end{eqnarray}
where $C=G\otimes\I_{p},~\mu=(\mu^{T}_{1},\ldots,\mu^{T}_{k})^{T}$, where $\otimes$ denotes the Kronecker product operator and $I_{p}$ denotes the p$\times$p identity matrix. For any $\delta\in\R^{p}$, (\ref{H:2.2}) can be transformed as follows:
\begin{align*}
C\mu=0\Leftrightarrow(I_{q}\otimes\delta^{T})C\mu=0.
 \end{align*}
Thus, the statistic for test (\ref{H:2.2}) based on the L$^{2}$-norm construction is:
\begin{align}\label{eq2.2}
T_{n}=\|(I_{q}\otimes\delta^{T})C\hat{\mu}\|^{2},
\end{align}
where $\hat{\mu}=(\bar{y}_{1}^{T},\ldots,\bar{y}_{k}^{T})^{T}$ and $\bar{y_{i}}$ is an unbiased estimate of $\mu_{i}$, $i\in\{1,\ldots\,k\}$ and $I_{q}$ denotes the q$\times$q identity matrix.similar to \cite{Jiang:2022},
\begin{align*}
T_{n}=\|(I_{q}\otimes\delta^{T})C\hat{\mu}\|^{2}&\Leftrightarrow(C\hat{\mu})^{T}(I_{q}\otimes\delta\delta^{T})C\hat{\mu}\\
&\Leftrightarrow\int(C\hat{\mu})^{T}(I_{q}\otimes\delta\delta^{T})(C\hat{\mu})w(\delta)d\delta.
 \end{align*}
where $w(\delta)$ denotes a positive weight. This transformation is essential for constructing the new test statistic, and we can obtain an explicit expression for $T_{n}$, which is given in the following theorem.
\begin{theorem}\label{th2.1}
Let $w(\delta)=\prod\limits_{i=1}^{p}w_{i}(\delta_{i})$, and $w_{i}(\cdot)$ denotes a density function with mean $\alpha_{i}$ and variance $\beta_{i}^{2}$, $i\in\{1,\ldots\,k\}$, then we have
\begin{align}\label{eq2.3}
T_{n}=\|(I_{q}\otimes\delta^{T})C\hat{\mu}\|^{2}=\hat{\mu}^{T}(H\otimes\W_{p})\hat{\mu},
\end{align}
where $H:(h_{ij})_{i,j=1}^{k}=G^{T}G=\tilde{G}^{T}(\tilde{G}D\tilde{G}^{T})^{-1}\tilde{G}$ and $W_{p}=B+aa^{T}$, where $a=(\alpha_{1},\ldots,\alpha_{p})^{T}$,

$\\$
$B=\begin{pmatrix}
\beta_{1}^{2} & 0 & \cdots & 0 \\
0 & \beta_{2}^{2} & \cdots & 0 \\
\vdots & \vdots & \ddots & \vdots \\
0 & 0 & \cdots & \beta_{p}^{2}
\end{pmatrix}$.
\end{theorem}
\begin{remark}
Theorem \ref{th2.1} is crucial in that it shows that to show that the null hypothesis holds is to show that $T_{n}=0$. When $\delta$ follows a density function with independent components, different tests can be produced by choosing different parameters for $\alpha_{i}$ and $\beta_{i}, i\in\{1,\ldots\,p\}$.
\end{remark}
To better study the statistic, we transform $T_{n}$ into
\begin{align}\label{eq2.4}
T_{n}=T_{n,0}+2S_{n}+\mu^{T}(H\otimes\W_{p})\mu,
\end{align}
where $T_{n,0}=(\hat{\mu}-\mu)^{T}(H\otimes\W_{p})(\hat{\mu}-\mu),~S_{n}=\mu^{T}(H\otimes\W_{p})(\hat{\mu}-\mu)$.
Under the null assumption that $T_{n,0}$ and $T_{n}$ have the same distribution, it follows from (\ref{eq2.1}) that H is held constant under this transformation, and hence $T_{n}$, $T_{n,0}$, $S_{n}$ and $\mu^{T}(H\otimes\W_{p})\mu$ are invariant.
For further study, we now set $\varphi=(\varphi_{1}^{T},\ldots,\varphi_{k}^{T})^{T}$ and let $\varphi_{i}=\sqrt{n_{i}}(\bar{y}_{i}-\mu_{i}),~i\in\{1,\ldots\,k\}$. By some simple calculations, we have $E(\varphi)=0_{kp},~Cov(\varphi)=\Sigma=diag(\Sigma_{1},\ldots,\Sigma_{k})_{kp\times\\kp}$. It's very simple to prove out that
\begin{align}\label{eq2.5}
T_{n,0}=\varphi^{T}(B^{T}B{\otimes}W_{p})\varphi=\varphi^{T}(A\otimes\W_{p})\varphi,
\end{align}
where $B=(\tilde{G}D\tilde{G}^{T})^{-1/2}\tilde{G}D^{1/2},~A:(a_{ij})_{i,j=1}^{k}=B^{T}B=D^{1/2}\tilde{G}^{T}(\tilde{G}D\tilde{G}^{T})^{-1}\tilde{G}D^{1/2}$. It's easy to see that $BB^{T}=I_{q}$ and A is an idempotent matrix with $A=A^{T},~A=A^{2}$ and $tr(A)=q$. Thus we have $a_{ii}>0,~i\in\{1,\ldots,k\}$.
\begin{theorem}\label{th2.2}
Let $\chi_{v}^{2}$ denote a central chi-square distribution and its degree of freedom is v. Set $A_{1},\ldots,A_{r},\ldots$ denote i.i.d. random variables. When all k samples are normal, we get $\varphi\backsim\N(0,\Sigma)$. For any given n and p, one obtains that the distribution of $T_{n,0}$ is the same as the following chi-square type mixtures
\begin{align}\label{eq2.6}
T_{n,0}^{*}=\sum_{r=1}^{kp}\lambda_{n,r}A_{r},
\end{align}
where $\lambda_{n,r},~r\in\{1,\ldots,(kp)\}$ is the descending eigenvalue of $\Omega_{n}=Cov{(B\otimes\W_{p})\varphi}=(B\otimes\W_{P})\Sigma(B^{T}\otimes\W_{P})$.
\end{theorem}
In order to better study $T_{n,0}^{*}$, by some algebraic calculations, we have
\begin{align}\label{eq2.7}
E(T_{n,0}^{*})=tr(\Omega_{n}),~Var(T_{n,0}^{*})=2tr(\Omega_{n}^{2}),~E{T_{n,0}^{*}-E(T_{n,0}^{*})}^{3}=8tr(\Omega_{n}^{3}).
\end{align}
Thus,the skewness of $T_{n,0}^{*}$ is given by
\begin{align}\label{eq2.8}
E\{T_{n,0}^{*}-E(T_{n,0}^{*})\}^{3}/Var^{3/2}(T_{n,0}^{*})=(8/d^{*})^{1/2},
\end{align}
where $d^{*}=tr^{3}(\Omega_{n}^{2})/tr^{2}(\Omega_{n}^{3})$. Observe that $T_{n,0}^{*}$ is often skewed and always nonnegative, though it can occasionally become asymptotically normal. Since $T_{n,0}^{*}$ is obtained from the test statistic $T_{n,0}$ when the k samples are normally distributed, we refer to the distribution of $T_{n,0}^{*}$ as the normal-reference distribution of $T_{n,0}$ accordance with \cite{Zhang:2021}. We can also demonstrate that, under certain regularity assumptions, $T_{n,0}$ and $T_{n,0}^{*}$ have the same normal or non-normal limit.

\subsection{Main results}\label{sec2.2}
\noindent In this subsection, the main conclusions of the paper are stated, and in order to obtain the properties of the test statistic, denote $\varrho_{n,r},~r\in\{1,\ldots,(kp)\}$ is the decreasing-ordered eigenvalues of $\Omega_{n}/\{tr(\Omega_{n}^{2})\}^{1/2}$, i.e., $\varrho_{n,r}=\lambda_{n}/\{tr(\Omega_{n}^{2})\}^{1/2},~r\in\{1,\ldots,(kp)\}$. For further study, we assume the following five conditions:
\begin{enumerate}
\item[(C1)] There exist p$\times$m matrix $\Gamma_{i}$ satisfies $\Gamma_{i}\Gamma_{i}^{T}=\Sigma_{i}$ and $z_{ij}'s$ are i.i.d. m-vectors, with $E(z_{ij})=0,~Cov(z_{ij})=I_{m}$, we let $y_{ij}=\mu_{i}+\Gamma_{i}z_{ij},~i\in\{1,\ldots,k\},~j\in\{1,\ldots,n_{i}\}$.

\item[(C2)] We assume $z_{ijl}$ denote the $l$-th component of $z_{ij}$. If there is one $v_{l}=1$ (two $v_{l}=2$) whenever $v_{1}+\cdots+v_{m}=4$, we have $E(z_{ijl}^{4})=3+\Delta<\infty$, and $E(z_{ij1}^{v1}{\ldots}z_{ijm}^{vm})=0~(or 1)$, where $\Delta$ is constant and $v_{1},\ldots,v_{m}$ are nonnegative integers.

\item[(C3)] As $n\rightarrow\infty$, $n_{i}/n\rightarrow\tau_{i}\in(0,1),~i\in\{1,\ldots,k\}$.

\item[(C4)] Assume $\lim_{n,p\rightarrow\infty}\varrho_{n,r}=\varrho_{r},~r\in\{1,2,\ldots\}$, and $\lim_{n,p\rightarrow\infty}\sum_{r=1}^{p}\varrho_{n,r}=\sum_{r=1}^{\infty}\varrho_{r}<\infty$.
    
\item[(C5)] For $i_{1},i_{2},i_{3},i_{4}\in\{1,\ldots,k\}$, as $p\rightarrow\infty$, $tr(W_{p}\Sigma_{i1}W_{p}\Sigma_{i2}W_{p}\Sigma_{i3}W_{p}\Sigma_{i4})=o\{tr(W_{p}\Sigma_{i1}\\W_{p}\Sigma_{i2})tr(W_{p}\Sigma_{i3}W_{p}\Sigma_{i4})\}$.
\end{enumerate}

\begin{remark}
As can be seen from the above conditions, we are not imposing a direct relationship between the data dimension p and the sample size n. Conditions C1 and C2 are actually extensions of those assumptions by \cite{Bai and Saranadasa:1996}. They are imposed to combat the non-normality of k high-dimensional samples. Condition C3 is a standard regularity assumption in k sample problems, which guarantees that k sample go to infinity proportionally. Condition C4 ensures that the limit of $\lambda_{n,r}$ exists when $n,p$ tends to infinity and that the limit and summation are interchangeable in expression $\lim_{n,p\rightarrow\infty}\sum_{r=1}^{p}\varrho_{n,r}$. It is used to ensure that the normalized $T_{n,0}$ and $T_{n,0}^{*}$ limit distributions are not normal. Condition C5 guarantees the consistency and asymptotic normality of our proposed test, ensuring that the limiting distributions of normalized $T_{n,0}$ and $T_{n,0}^{*}$ are normal.
 \end{remark}
Let $\overset{d}{=}$ denote  equality in distribution and $\overset{L}{\longrightarrow}$ denote convergence in distribution.
\begin{theorem}\label{th2.3}~\\
(1) Under Conditions C1-C4, when $n,p\rightarrow\infty$, we have
\begin{align}\label{eq2.9}
\frac{T_{n.0}-tr(\Omega_{n})}{\sqrt{2tr(\Omega_{n}^{2})}}\overset{L}{\longrightarrow}\zeta,~\frac{T_{n.0}^{*}-tr(\Omega_{n})}{\sqrt{2tr(\Omega_{n}^{2})}}\overset{L}{\longrightarrow}\zeta,
\end{align}
where $\zeta\overset{d}{=}\sum_{r=1}^{\infty}\varrho_{r}(A_{r}-1)/\sqrt{2}$.\\
(2) Under Conditions C1-C3 and C5, when $n,p\rightarrow\infty$, we have
\begin{align}\label{eq2.10}
\frac{T_{n.0}-tr(\Omega_{n})}{\sqrt{2tr(\Omega_{n}^{2})}}\overset{L}{\longrightarrow}N(0,1),~\frac{T_{n.0}^{*}-tr(\Omega_{n})}{\sqrt{2tr(\Omega_{n}^{2})}}\overset{L}{\longrightarrow}N(0,1).
\end{align}
Then under the conditions of (1) or (2), we always have
\begin{align}\label{eq2.11}
\sup_{x}|Pr(T_{n,0}\leq\x)-Pr(T_{n,0}^{*}\leq\x)|\longrightarrow0.
\end{align}
\end{theorem}
\begin{remark}
Theorem \ref{th2.3} provides a systematic justification for our use of the distribution of $T_{n,0}^{*}$ to approximate the distribution of $T_{n,0}$. It can be shown that the asymptotic distribution of $T_{n,0}$ depends in a complicated way on the limiting ratio of group sample size to total sample size and the group covariance matrix. When the k group covariance matrices are the same, the asymptotic distribution of $T_{n,0}$ depends only on the common covariance matrix.
 \end{remark}

Theorem \ref{th2.3} shows that, in order to realize the proposed test, we can approximate the distribution of $T_{n,0}$ by the distribution of $T_{n,0}^{*}$. From (\ref{eq2.6}) it is known that $T_{n,0}^{*}$ is a cartesian mixture and that the coefficients $\lambda_{n,r}$ are unknown and it's the eigenvalues of $\Omega_{n}$. $T_{n,0}^{*}$ is nonnegative and usually skewed. Therefore, it is not natural to always approximate its distribution by a normal distribution, as many methods mentioned in the literature do. Similar to what was mentioned in \cite{Zhang et al:2022(b)}, it is natural to approximate its distribution by the well-known Box $\chi^{2}$-approximation, also known as the Welch-Satterthwaite $\chi^{2}$-approximation. We approximate the distribution of $T_{n,0}^{*}$ by a random variable of the following form:
\begin{align}\label{eq2.12}
R\overset{d}{=}\beta\chi_{d}^{2}.
\end{align}
The parameters $\beta$ and $d$ are determined by matching the first two cumulants, i.e., the mean and variance, of $T_{n,0}^{*}$ and $R$. The first two cumulants of $T_{n,0}^{*}$ are in (\ref{eq2.7}), and the first two cumulants of $R$ are $E(R)={\beta}d$ and $Var(R)=2\beta^{2}d$. Equalizing the first two cumulants of $T_{n,0}^{*}$ and $R$ yields the following equation for $\beta$ and $d$:
\begin{align}\label{eq2.13}
\beta=\frac{tr(\Omega_{n}^{2})}{tr(\Omega_{n})},~d=\frac{tr^{2}(\Omega_{n})}{tr(\Omega_{n}^{2})},
\end{align}
where $\beta$ and $d$ are approximate parameters and $d$ is its approximate degree of freedom.The W-S $\chi^{2}$-approximation is very accurate and has been widely adopted as an approximate solution to many heteroskedastic problems in the classical setting. There are several advantages to applying the WS $\chi^{2}$-approximation: first, it is simple to implement and very fast to compute. The details of the computation will be mentioned below, and it is relatively simple and fast as it only requires the computation of some simple forms of estimators. Second, the W-S $\chi^{2}$-approximation guarantees that $T_{n,0}^{*}$ and $R$ have the same mean, variance, range, and similar shape. On the contrary, the normal approximation only guarantees the same mean and variance. And the proposed $L^{2}$-norm based W-S $\chi^{2}$-approximation test is expected to outperform the existing competitors under normal approximation in terms of size control. Last but not least, the degrees of freedom $d$ of the W-S $\chi^{2}$-approximation adapt to the shape of the distribution of $T_{n,0}^{*}$. According to Theorem 4 mentioned in \cite{Zhang:2020}, we know that $T_{n,0}^{*}$ is asymptotically normal when $d^{*}\rightarrow\infty$. Also, according to Theorem 5 in \cite{Zhang:2020}, we have $1{\leq}d^{*}{\leq}d{\leq}p$. Thus, when $T_{n,0}^{*}$ is asymptotically normal, at that time $d,d^{*}\rightarrow\infty$ and $R$ is also asymptotically normal, and when $d$ is asymptotically bounded, so is $d^{*}$, and thus neither $T_{n,0}^{*}$ nor $R$ is asymptotically normal. So W-S $\chi^{2}$-approximation is better adaptive than the normal approximation.
\begin{theorem}\label{th2.4}
Under Conditions C1-C3 and C5, when $n,p\rightarrow\infty$, we have $d\rightarrow\infty$, and
\begin{align}\label{eq2.14}
\frac{R-tr(\Omega_{n})}{\sqrt{2tr(\Omega_{n}^{2})}}\overset{L}{\longrightarrow}N(0,1),
\end{align}
\begin{align}\label{eq2.15}
\sup_{x}|Pr(T_{n,0}\leq\x)-Pr(R\leq\x)|\longrightarrow0.
\end{align}
\end{theorem}
\begin{remark}
Theorem \ref{th2.3} and the above theorem show that under Conditions C1-C3 and C5, $T_{n,0}$, $T_{n,0}^{*}$, and $R$ are asymptotically normal, and the W-S $\chi^{2}$-approximation is equal to the normal approximation. However, under Conditions C1-C4, Theorem \ref{th2.3} shows that $T_{n,0}$ and $T_{n,0}^{*}$ are asymptotically skewed, and thus the W-S $\chi^{2}$-approximation should be preferred to the normal approximation in this case.
\end{remark}

To formulate a test procedure, we need consistent estimates for $tr(\Omega_{n}),~tr^{2}(\Omega_{n})$ and $tr(\Omega_{n}^{2})$. By using the expression for $\Omega_{n}$ above, we get
\begin{equation}
\begin{aligned}
&tr(\Omega_{n})=\sum_{i=1}^{k}a_{ii}tr(W_{p}\Sigma_{i}),\\
&tr^{2}(\Omega_{n})=\sum_{i=1}^{k}a_{ii}^{2}tr^{2}(W_{p}\Sigma_{i})+2\sum_{1{\leq}i<j{\leq}k}a_{ii}a_{jj}tr(W_{p}\Sigma_{i})tr(W_{p}\Sigma_{j}),\\
&tr(\Omega_{n}^{2})=\sum_{i=1}^{k}a_{ii}^{2}tr\{(W_{p}\Sigma_{i})^{2}\}+2\sum_{1{\leq}i<j{\leq}k}a_{ii}^{2}tr(W_{p}\Sigma_{i}W_{p}\Sigma_{j}),
\end{aligned}
\label{eq2.16}
\end{equation}
\allowdisplaybreaks
where $a_{ij},~i,j\in\{1,\ldots,k\}$ denotes the entries of A. It involves all the group covariance matrices $\Sigma_{i}$. In order to estimate $tr(\Omega_{n}),tr^{2}(\Omega_{n})$ and $tr(\Omega_{n}^{2})$ consistently, we must first estimate $tr(W_{p}\Sigma_{i}),~tr(W_{p}\Sigma_{i})tr(W_{p}\Sigma_{j}),~tr(W_{p}\Sigma_{i}W_{p}\Sigma_{j})$ and $tr^{2}(W_{p}\Sigma_{i}),~tr\{(W_{p}\Sigma_{i})^{2}\}$ consistently. We denote the consistent estimates of $tr(W_{p}\Sigma_{i}),~tr(W_{p}\Sigma_{i})tr(W_{p}\Sigma_{j}),~tr(W_{p}\Sigma_{i}W_{p}\Sigma_{j})$ by $tr(W_{p}\widehat{\Sigma}_{i}),~tr(W_{p}\widehat{\Sigma}_{i})tr(W_{p}\widehat{\Sigma}_{j}),~tr(W_{p}\widehat{\Sigma}_{i}W_{p}\widehat{\Sigma}_{j})$, respectively, where $\widehat{\Sigma}_{i}=(n_{i}-1)^{-1}\sum_{i=1}^{n_{i}}(\\y_{ij}-\bar{y_{i}})(y_{ij}-\bar{y_{i}})^{T}$ is the usual unbiased estimator of $\Sigma_{i},~i\in\{1,\ldots,k\}$, respectively. Carefully studying equation (27) of \cite{Zhang et al:2022(b)}, Similarly, we obtain the ratio-consistent estimators of $tr^{2}(W_{p}\Sigma_{i})$ and $tr\{(W_{p}\Sigma_{i})^{2}\}$,
\begin{equation}
\begin{aligned}
&\Widehat{tr^{2}(W_{p}\Sigma_{i})}=\frac{(n_{i}-1)n_{i}}{(n_{i}-2)(n_{i}+1)}\left\{tr^{2}(W_{p}\widehat{\Sigma}_{i})-\frac{2}{n_{i}}tr\{(W_{p}\widehat{\Sigma}_{i})^{2}\}\right\},\\
&\Widehat{tr\{(W_{p}\Sigma_{i})^{2}\}}=\frac{(n_{i}-1)^{2}}{(n_{i}-2)(n_{i}+1)}\left\{tr\{(W_{p}\widehat{\Sigma}_{i})^{2}\}-\frac{1}{n_{i}-1}tr^{2}(W_{p}\widehat{\Sigma}_{i})\right\},
\end{aligned}
\label{eq2.17}
\end{equation}
where $i,j\in\{1,\ldots,k\}$. Thus, we obtain the ratio-consistent estimators of $tr(\Omega_{n}),~tr^{2}(\Omega_{n})$ and $tr(\Omega_{n}^{2})$,
\begin{equation}
\begin{aligned}
&tr(\widehat{\Omega}_{n})=\sum_{i=1}^{k}a_{ii}tr(W_{p}\widehat{\Sigma_{i}}),\\
&\widehat{tr^{2}(\Omega_{n})}=\sum_{i=1}^{k}a_{ii}^{2}\Widehat{tr^{2}(W_{p}\Sigma_{i})}+2\sum_{1{\leq}i<j{\leq}k}a_{ii}a_{jj}tr(W_{p}\widehat{\Sigma_{i}})tr(W_{p}\widehat{\Sigma_{j}}),\\
&\widehat{tr(\Omega_{n}^{2})}=\sum_{i=1}^{k}a_{ii}^{2}\Widehat{tr\{(W_{p}\Sigma_{i})^{2}\}}+2\sum_{1{\leq}i<j{\leq}k}a_{ii}^{2}tr(W_{p}\widehat{\Sigma_{i}}W_{p}\widehat{\Sigma_{j}}).
\end{aligned}
\label{eq2.18}
\end{equation}
Therefore, we obtain the estimators of $\beta$ and $R$,
\begin{align}\label{eq2.19}
\widehat{\beta}=\frac{\widehat{tr(\Omega_{n}^{2})}}{tr(\widehat{\Omega}_{n})},~\widehat{d}=\frac{\widehat{tr^{2}(\Omega_{n})}}{\widehat{tr(\Omega_{n}^{2})}}.
\end{align}
Let $\chi_{d}^{2}(\alpha)$ denote the upper $100\alpha$ percentile of $\chi_{d}^{2}$, where any nominal significance level $\alpha>0$. Therefore, we have the following theorem.
\begin{theorem}\label{th2.5}
Under Conditions C1-C3, when $n\rightarrow\infty$, we have
\begin{align}\label{eq2.20}
\frac{tr(\widehat{\Omega}_{n})}{tr(\Omega_{n})}\overset{P}{\longrightarrow}1,~\frac{\widehat{tr^{2}(\Omega_{n})}}{tr^{2}(\Omega_{n})}\overset{P}{\longrightarrow}1,~\frac{\widehat{tr(\Omega_{n}^{2})}}{tr(\Omega_{n}^{2})}\overset{P}{\longrightarrow}1,
\end{align}
and\\
\begin{align}\label{eq2.21}
\frac{\widehat{\beta}}{\beta}\overset{P}{\longrightarrow}1,~\frac{\widehat{d}}{d}\overset{P}{\longrightarrow}1,~\frac{\widehat{\beta}\chi_{\widehat{d}}^{2}(\alpha)}{\beta\chi_{d}^{2}(\alpha)}\overset{P}{\longrightarrow}1,
\end{align}
uniformly for all p.
\end{theorem}
\begin{remark}
Theorem \ref{th2.5} shows that the conclusion in Theorem \ref{th2.3} still holds when $tr(\Omega_{n}),~tr^{2}(\Omega_{n})$ and $tr(\Omega_{n}^{2})$ is replaced by its ratio-consistent estimators.
\end{remark}
\subsection{Power of the proposed test}\label{sec2.3}
\noindent In this subsection, we investigate the power of the proposed test. By some simple algebraic calculations, we have
\begin{align}\label{eq2.22}
E(S_{n})=0,~Var(S_{n})=\mu^{T}(H{\otimes}W_{p})diag(\frac{\Sigma_{1}}{n_{1}},\ldots,\frac{\Sigma_{k}}{n_{k}})(H{\otimes}W_{p})\mu.
\end{align}
For simplicity, we consider the asymptotic powers of $T_{n}$ under the following local alternatives. When $n,p\rightarrow\infty$, we have
 \begin{align}\label{eq2.23}
\frac{Var(S_{n})}{Var(T_{n,0})}=o(1).
\end{align}
Under this local alternatives, we have $\frac{S_{n}}{\sqrt{Var(T_{n,0})}}\overset{P}{\longrightarrow}0.$ Then under Condition C3, when $n\rightarrow\infty$, we have
 \begin{align}\label{eq2.24}
\Omega_{n}\longrightarrow\Omega=(B^{*}{\otimes}W_{p})\Sigma(B^{*T}{\otimes}W_{p}),
\end{align}
where $\lim_{n\rightarrow\infty}B=B^{*}$, and denotes $H^{*}=\lim_{n\rightarrow\infty}n^{-1}H$.The following theorem gives the asymptotic power function of $T_{n}$.
\begin{theorem}\label{th2.6}~\\
(1) Under Conditions C1-C4 and the local alternative (\ref{eq2.19}), when $n,p\rightarrow\infty$, we have
\begin{align}\label{eq2.25}
Pr\left\{T_{n}>\widehat{\beta}\chi_{\widehat{d}}^{2}(\alpha)\right\}=Pr\left\{\zeta\geq\frac{\chi_{d}^{2}(\alpha)-d}{\sqrt{2d}}-\frac{ntr(W_{p}M^{T}H^{*}M)}{\sqrt{2tr(\Omega^{2})}}\right\}\{1+o(1)\},
\end{align}
where $\zeta$ is defined in Theorem \ref{th2.3}(1).\\
(2) Under Conditions C1-C3, C5 and the local alternative (\ref{eq2.19}), when $n,p\rightarrow\infty$, we have
\begin{align}\label{eq2.26}
Pr\left\{T_{n}>\widehat{\beta}\chi_{\widehat{d}}^{2}(\alpha)\right\}=\Phi\left\{-z_{\alpha}+\frac{ntr(W_{p}M^{T}H^{*}M)}{\sqrt{2tr(\Omega^{2})}}\right\}\{1+o(1)\},
\end{align}
where $z_{\alpha}$ denotes the upper $100\alpha$-percentile of $N(0,1)$.
\end{theorem}

\section{Simulation Studies}\label{sec3}
\noindent In this section,  we present four simulation studies. Simulation \ref{sec3.1} is designed to compare the empirical powers of $T_{n}$ and existing tests for the MANOVA problem in different situations. Simulation \ref{sec3.2} is designed to demonstrate how $T_{n}$ and existing tests perform the GLHT problem in different situations, as well as a comparison of the empirical powers of $T_{n}$ and existing tests. Similar to what was mentioned in \cite{Jiang:2022}, we set $\alpha_{1}=\cdots=\alpha_{p}=2p^{-3/8}$ and $\beta_{i}=\sqrt{2}(p+i)/p,~i = \{1,\ldots,p\}$ for convenience. In both simulations, we set the nominal significance level to $\alpha$=0.05. In this section, denote the test of the normal reference $L^{2}$-norm based on the W-S $\chi^{2}$-approximation studied in Section \ref{sec2} as $T_{n}$. For simplicity, the k high-dimensional samples are generated from the following model,
\begin{align}\label{eq3.1}
\mathbf{y}_{ij}=\mu_{i}+\Gamma_{i}\mathbf{z}_{ij},~i\in\{1,\ldots,\},~j\in\{1,\ldots,\}
\end{align}
where $\mathbf{z}_{ij}=(\mathbf{z}_{ij1},\ldots,\mathbf{z}_{ijp})^{T},~i\in\{1,\ldots,\},~j\in\{1,\ldots,\}$ are independent m-dimensional random variables. To better study the power of $T_{n}$, we set $\mathbf{z}_{ij}$ to arise from the following three models, representing normal, non-normal but symmetric, and non-normal and skewed distribution shapes, respectively,
\begin{flalign*}
~~~~&Model~1:~\mathbf{z}_{ijl}\overset{i.i.d}{\sim}N(0,1),~l\in\{1,\ldots,p\}.& \\
&Model~2:~\mathbf{z}_{ijl}=\omega_{ijl}/\sqrt{2},~\omega_{ijl}\overset{i.i.d}{\sim}t_{4},~l\in\{1,\ldots,p\}.& \\
&Model~3:~\mathbf{z}_{ijl}=(\omega_{ijl}-1)/\sqrt{2},~\omega_{ijl}\overset{i.i.d}{\sim}\chi_{1}^{2},~l\in\{1,\ldots,p\}.&
\end{flalign*}
In both simulations, we set three sets of group sample sizes: $n_{1}=(n_{1},n_{2},n_{3},n_{4})=(25,30,40,45),n_{2}=(30,40,55,60),n_{3}=(80,95,115,120)$, and three sets of dimensions: $p=(100,200,300,400,500)$. We assigned $\mu_{1}=\mu_{2}=\mu_{3}=\mu_{4}=0$ in the case of $H_{0}$, $\mu_{1}=\mu_{2}=\mu_{3}=0$,~$\mu_{4}$ had $p^{1-\rho}$ non-zero entries of equal value that were uniformly allocated among $\{1,\ldots,p\}$ in the case of $H_{1}$, and the values of the nonzero entries were $\sqrt{2r(1/n_{1}+1/n_{2}+1/n_{3}+1/n_{}4){\log}~p}$, where $r=(0.02,0.04,0.06,0.08)$ controls the signal strength, $\rho=(0.1,0.2,0.3,0.4)$ controls the signal sparsity and covers highly dense signals when $\rho=0.1$, moderately dense signals when $\rho=0.2$ or $\rho=0.3$, and  moderately sparse when $\rho=0.4$. For the covariance matrix, we consider the following two scenarios:
\begin{description}
\item[\bf{Scenario 1:}] $\mathbf{\Sigma_{1}=(0.4^{|i-j|}),~\Sigma_{2}=I_{p},\Sigma_{3}=3I_{p},\Sigma_{4}=6I_{p}}$ for $1{\leq}i, j{\leq}p$ and $I_{p}$ denotes the $p{\times}p$ identity matrix. 
\item[\bf{Scenario 2:}] $\mathbf{\Sigma_{1}=3I_{p},\Sigma_{2}=5I_{p},\Sigma_{3}=0.09I_{p}+0.01J_{p},\Sigma_{4}=2\Sigma_{1}}$, where $I_{p}$ denotes the $p{\times}p$ identity matrix, $J_{p}$ denotes the matrix of ones of size $p{\times}p$.
\end{description}
The remaining cases are shown in the following two simulations.
\subsection{A simulation study of the MANOVA problem}\label{sec3.1}
\noindent In this simulation, we compete the new test statistic $T_{n}$ with the three existing tests on the MANOVA problem. We denote the test statistic proposed by \cite{Zhou et al:2017} as $T_{zb}$, the one proposed by \cite{Zhang et al:2022(b)} as $T_{zjt}$, and the one proposed by \cite{Hu et al:2017} as $T_{hj}$. In this paper, the hypothesis problem in \ref{H:1.2} reduces to the MANOVA problem when we set $\tilde{G}=(-1_{k-1},I_{k-1})$ or $\tilde{G}=(I_{k-1},-1_{k-1})$, where $1_{k-1}$ denotes the $(k-1)$-dimensional vector with entries 1 and $I_{k-1}$ denotes the $(k-1){\times}(k-1)$ identity matrix. For simplicity, we set $k=4$. Tables \ref{tab1} and \ref{tab2} demonstrate the empirical sizes in scenario 1 and scenario 2.

\noindent For each setting, 2000 replications are simulated to calculate all empirical sizes and power levels. Table \ref{tab1} and Table \ref{tab2} demonstrate the empirical sizes of the four test statistics for different scenarios. The nominal size is 0.05, as can be seen from Tables \ref{tab1} and Table \ref{tab2}, $T_{n}$ performs well in various cases in terms of size control, which suggests that our new test using the W-S $\chi^{2}$-approximation is working well. In addition, we find that the empirical sizes of the proposed test statistics are well controlled around 0.05.\\
The empirical powers of our proposed test statistic and the other three comparison statistics are displayed in Figures 1-6. The powers of $T_{n}$ are similar under the three different distributions, so we report only the empirical powers under $N(0,1)$ in Figures 1-6. For the other two distributions, the results are shown in the Supplementary Material. From Figures 1-6, we can derive that the proposed test is most powerful at $\rho=0.1$ or $\rho=0.2$. Thus, $T_{n}$ is powerful when nonzero signals of the difference between two mean vectors are weakly dense with nearly the same sign, and the empirical power of $T_{n}$ diminishes as $\rho$ increases. Besides, the empirical power of $T_{n}$ increases as the signal strength r increases; the empirical powers of all four tests were lower when $\rho=0.3$ or $\rho=0.4$.

\begin{table}[H]
\begin{center}
\caption{{Empirical sizes in the MANOVA problem when $z_{ij}$ follows $N(0,1)$ for scenario 1.}}
\label{tab1}
\begin{tabularx}{\textwidth}{XXXXXX}
\hline
$p$  &$n$  &\multicolumn{4}{c}{$\mathbf{z}_{ijl}\overset{i.i.d}{\sim}N(0,1)$} \\ 
\cmidrule(r){3-6}
                    &  &$T_{n}$ &$T_{zb}$  &$T_{zjt}$ &$T_{hj}$  \\
\hline
100    &$n_{1}$   &0.049	&0.059	&0.0535	&0.0565	\\
       &$n_{2}$   &0.049	&0.0585	&0.055	&0.059	 \\
       &$n_{3}$   &0.051	&0.0525	&0.0535	&0.054	 \\
200    &$n_{1}$   &0.0495	&0.049	&0.0455	&0.048	 \\
       &$n_{2}$   &0.0505	&0.044	&0.0465	&0.048	 \\
       &$n_{3}$   &0.053	&0.0565	&0.055	&0.058	 \\
300    &$n_{1}$   &0.052	&0.056	&0.055	&0.0565	 \\
       &$n_{2}$   &0.0525	&0.0585	&0.0585	&0.057	 \\
       &$n_{3}$   &0.0485	&0.043	&0.043	&0.0445	 \\
400    &$n_{1}$   &0.0505	&0.057	&0.0515	&0.055	 \\
       &$n_{2}$   &0.0515	&0.056	&0.051	&0.061	 \\
       &$n_{3}$   &0.053	&0.056	&0.053	&0.0565	 \\
500    &$n_{1}$   &0.0505	&0.047	&0.0455	&0.0525	 \\
       &$n_{2}$   &0.049	&0.051	&0.0485	&0.0545	 \\
       &$n_{3}$   &0.0515	&0.0615	&0.055	&0.0585	 \\
\hline
\end{tabularx}
\end{center}
\end{table}

 \begin{figure}[H]
 \centering
     \hspace{-0.8cm}
     \includegraphics[scale=1]{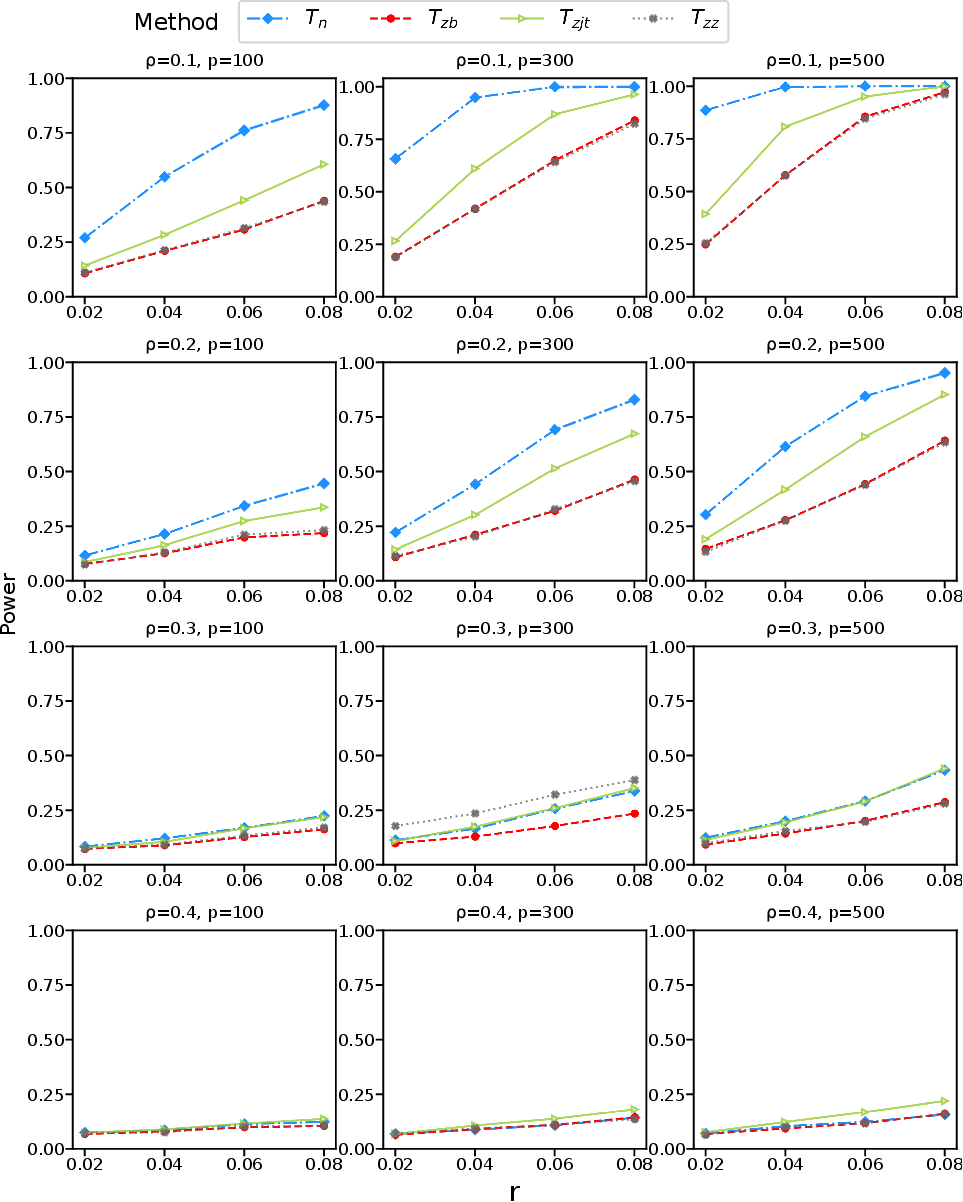}
    \begin{flushleft}
    {\par{\footnotesize\textbf{Figure 1.} Empirical powers when $z_{ij}$ follows $N(0,1)$ and $n_{1}=(25,30,40,45)$ for scenario 1 under different signal levels of $r$ and sparsity levels of $\rho$ in MANOVA problem.}}
     \end{flushleft}
 \end{figure}
 
 \begin{figure}[H]
 \centering
 \hspace{-0.8cm}
    \includegraphics[scale=1]{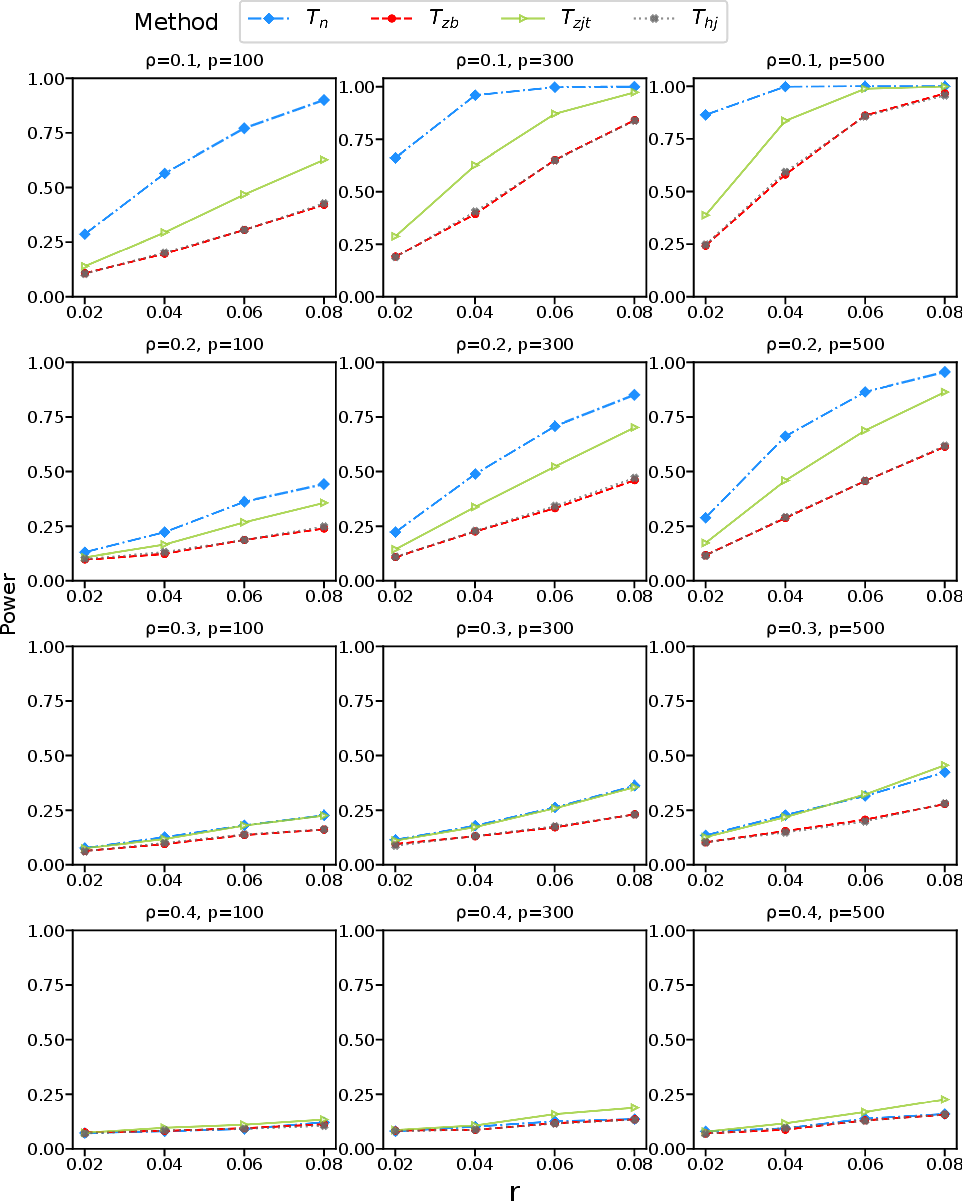}
      \begin{flushleft}
    {\par{\footnotesize\textbf{Figure 2.} Empirical powers when $z_{ij}$ follows $N(0,1)$ and $n_{2}=(30,40,55,60)$ for scenario 1 under different signal levels of $r$ and sparsity levels of $\rho$ in MANOVA problem.}}
 \end{flushleft}
 \end{figure}
 
  \begin{figure}[H]
  \centering
  \hspace{-0.8cm}
     \includegraphics[scale=1]{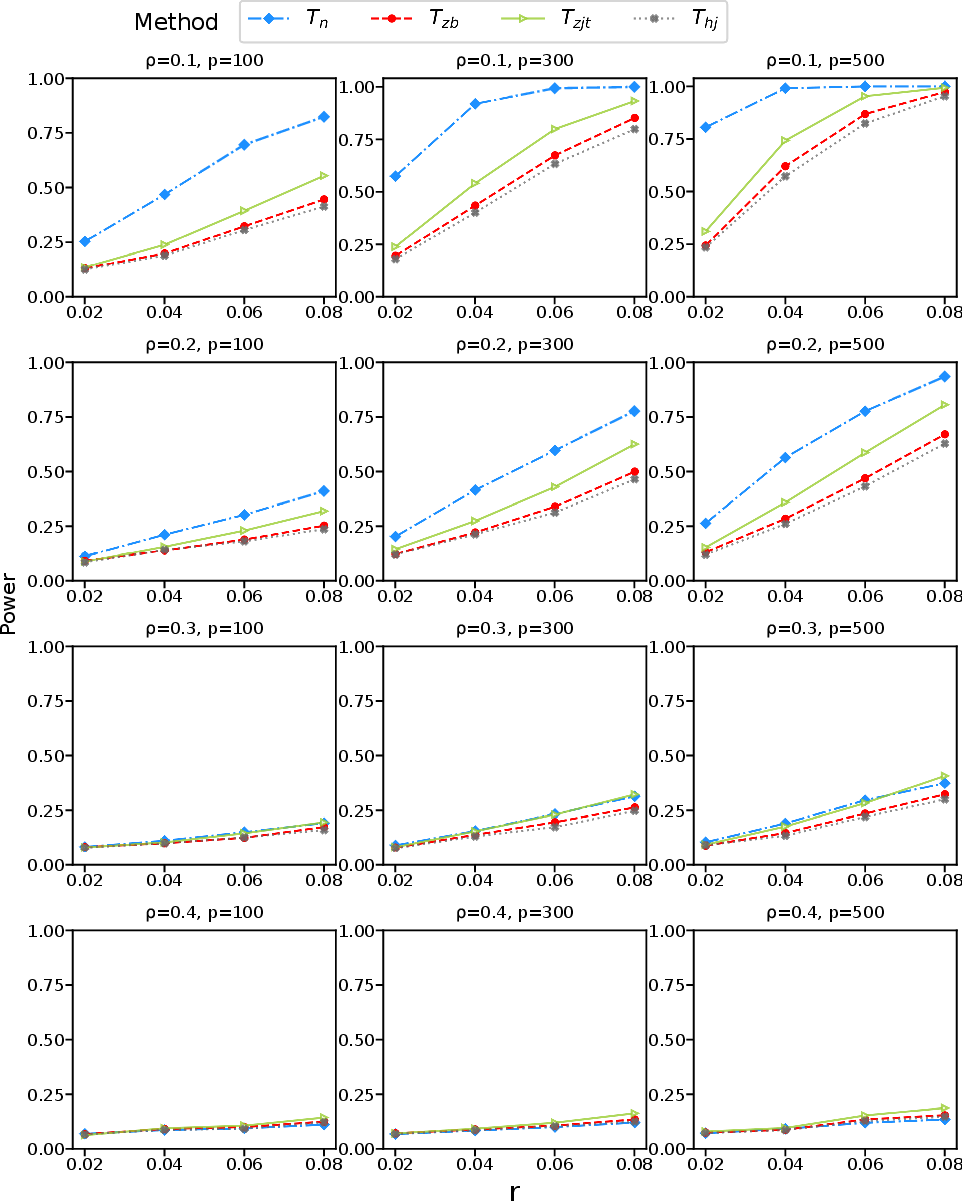}
     \begin{flushleft}
    {\par{\footnotesize\textbf{Figure 3.} Empirical powers when $z_{ij}$ follows $N(0,1)$ and $n_{3}=(80,95,115,120)$ for scenario 1 under different signal levels of $r$ and sparsity levels of $\rho$ in MANOVA problem.}}
     \end{flushleft}
 \end{figure}

\begin{table}[H]
\begin{center}
\caption{{Empirical sizes in the MANOVA problem when $z_{ij}$ follows $N(0,1)$ for scenario 2.}}
\label{tab2}
\begin{tabularx}{\textwidth}{XXXXXX}
\hline
$p$  &$n$  &\multicolumn{4}{c}{$\mathbf{z}_{ijl}\overset{i.i.d}{\sim}N(0,1)$} \\ 
\cmidrule(r){3-6} 
                    &  &$T_{n}$ &$T_{zb}$  &$T_{zjt}$ &$T_{hj}$  \\
\hline
100    &$n_{1}$   &0.0495	&0.049	&0.048	&0.0495			  \\
       &$n_{2}$   &0.053	&0.06	  &0.0585	&0.064			  \\
       &$n_{3}$   &0.0495	&0.06	  &0.056	&0.061			 \\
200    &$n_{1}$   &0.0485	&0.059	&0.055	&0.059			 \\
       &$n_{2}$   &0.049	&0.0555	&0.052	&0.063			   \\
       &$n_{3}$   &0.0525	&0.06	  &0.053	&0.057			  \\
300    &$n_{1}$   &0.054	&0.0565	&0.054	&0.0595			  \\
       &$n_{2}$   &0.053	&0.058	&0.056	&0.0585			  \\
       &$n_{3}$   &0.0525	&0.045	&0.0445	&0.0395				 \\
400    &$n_{1}$   &0.0515	&0.049	&0.048	&0.044			 \\
       &$n_{2}$   &0.0515	&0.048	&0.056	&0.054			\\
       &$n_{3}$   &0.0505	&0.049	&0.0505	&0.05		     \\
500    &$n_{1}$   &0.051	&0.0505	&0.0525	&0.0505			 \\
       &$n_{2}$   &0.0495	&0.061	&0.058	&0.058			  \\
       &$n_{3}$   &0.0495	&0.0425	&0.045	&0.053		 \\                                                                                                                 
\hline
\end{tabularx}
\end{center}
\end{table}

  \begin{figure}[H]
   \centering
   \hspace{-0.8cm}
    \includegraphics[scale=1]{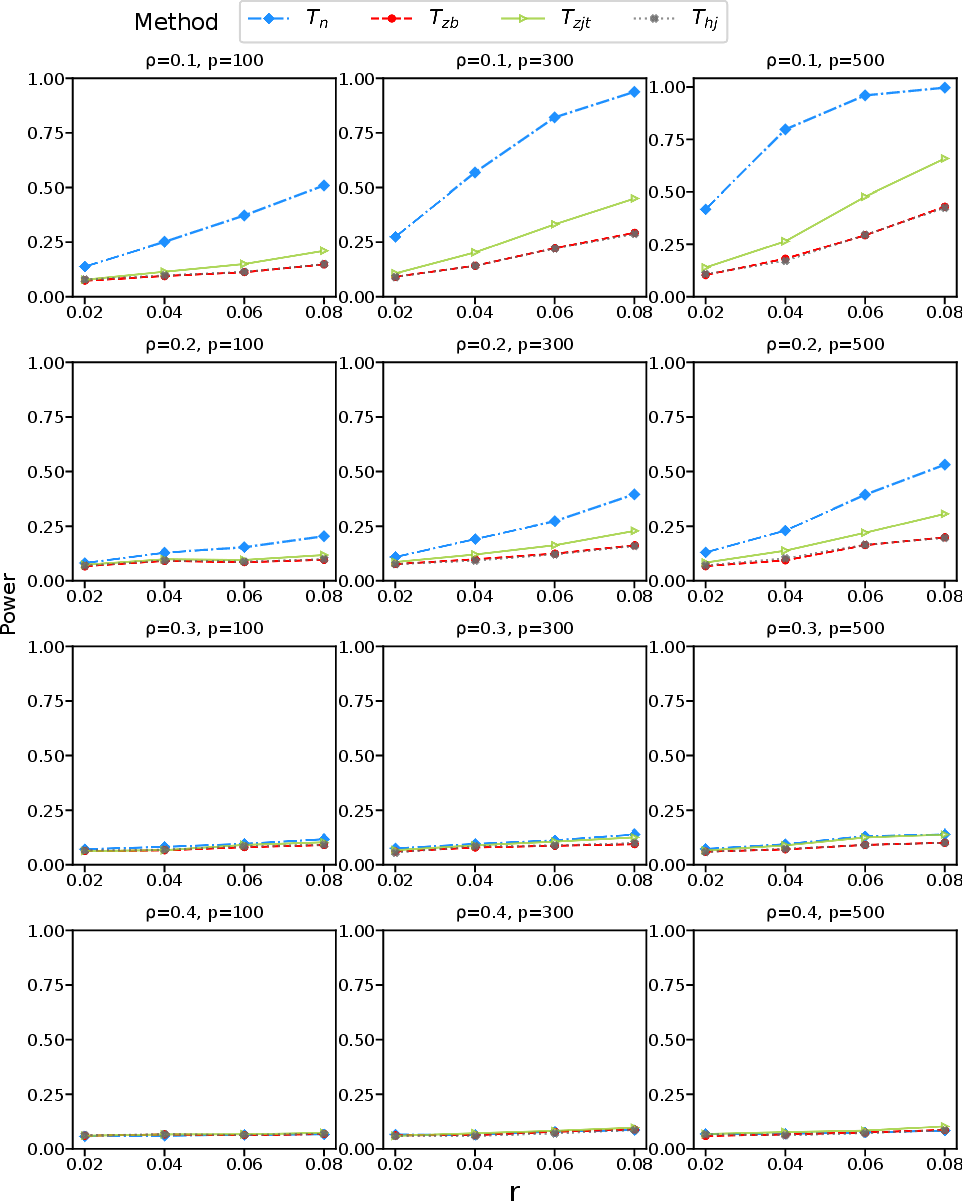}
    \begin{flushleft}
    {\par{\footnotesize\textbf{Figure 4.} Empirical powers when $z_{ij}$ follows $N(0,1)$ and $n_{1}=(25,30,40,45)$ for scenario 2 under different signal levels of $r$ and sparsity levels of $\rho$ in MANOVA problem.}}
     \end{flushleft}
 \end{figure}
 
  \begin{figure}[H]
   \centering
   \hspace{-0.8cm}
    \includegraphics[scale=1]{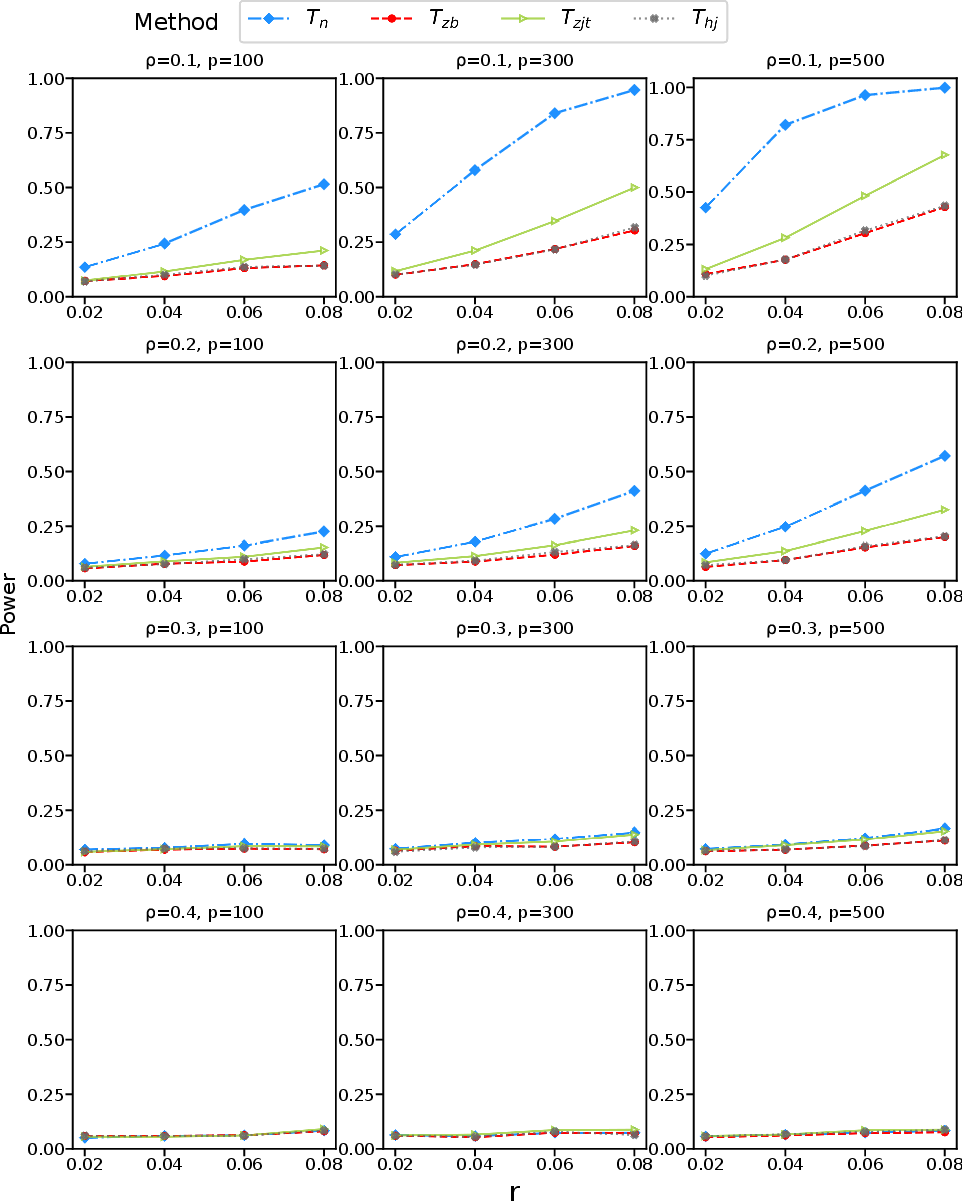}
    \begin{flushleft}
    {\par{\footnotesize\textbf{Figure 5.} Empirical powers when $z_{ij}$ follows $N(0,1)$ and $n_{2}=(30,40,55,60)$ for scenario 2 under different signal levels of $r$ and sparsity levels of $\rho$ in MANOVA problem.}}
     \end{flushleft}
 \end{figure}
 
  \begin{figure}[H]
  \centering
  \hspace{-0.8cm}
  \includegraphics[scale=1]{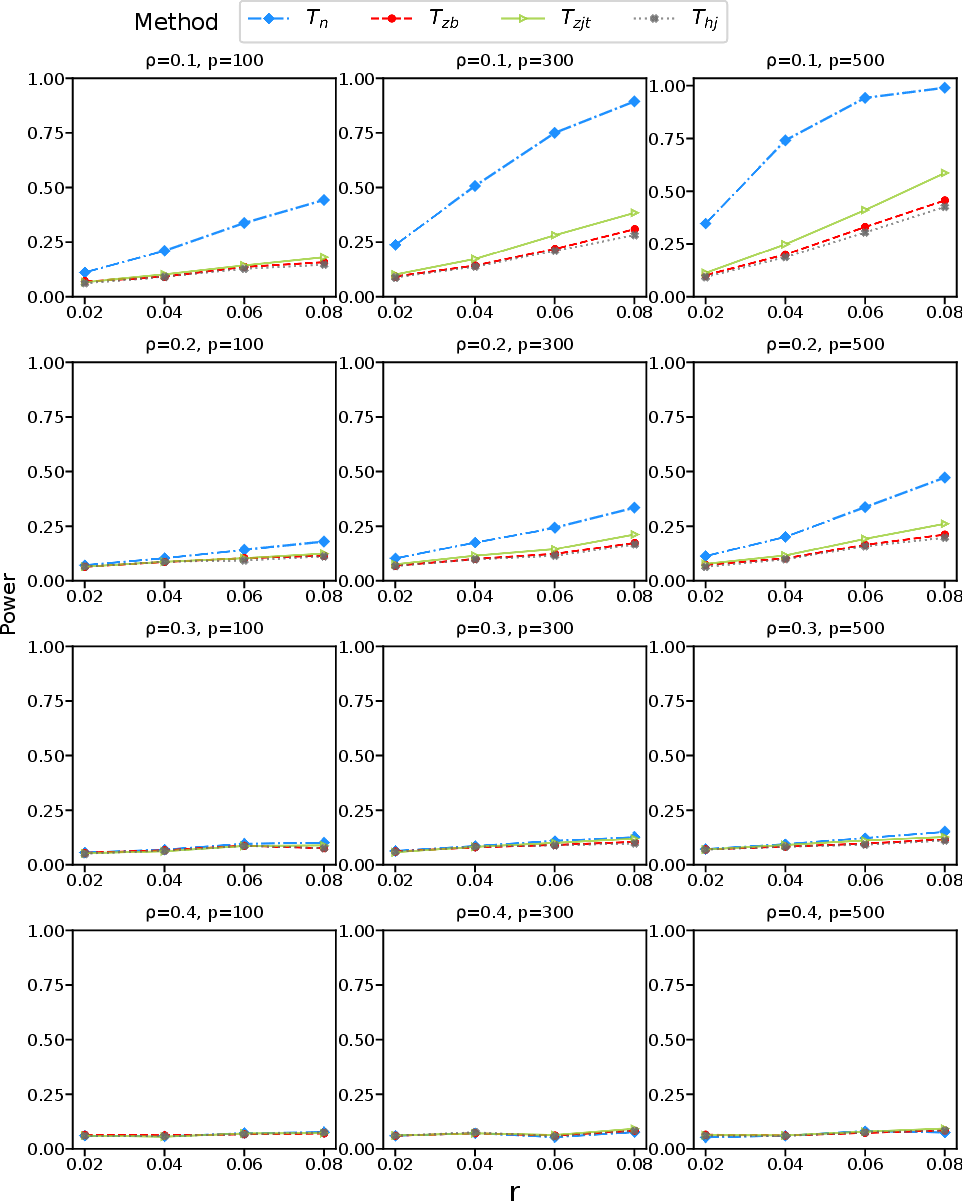}
  \begin{flushleft}
    {\par{\footnotesize\textbf{Figure 6.} Empirical powers of when $z_{ij}$ follows $N(0,1)$ and $n_{3}=(80,95,115,120)$ for scenario 2 under different signal levels of $r$ and sparsity levels of $\rho$ in MANOVA problem.}}
     \end{flushleft}
 \end{figure}

\subsection{A simulation study of the GLHT problem}\label{sec3.2}
\noindent In this simulation, we compete the new test statistic $T_{n}$ with the three existing tests on the GLHT problem. We denote the test statistic proposed by \cite{Zhou et al:2017} as $T_{zb}$, the one proposed by \cite{Zhang et al:2022(b)} as $T_{zjt}$, and the one proposed by \cite{Zhang et al:2022(a)} as $T_{zz}$. For convenience, we set $\tilde{G}=(e_{1,4}+2e_{2,4}+e_{3,4}-4e_{4,4})$. Tables 3 and 4 demonstrate the empirical sizes in
scenario 1 and scenario 2.

\noindent Similarly, for each setting, 5000 replications are simulated to calculate all empirical sizes and power levels. Table \ref{tab3} and Table \ref{tab4} demonstrate the empirical sizes of the four test statistics for different scenarios. The nominal size is also set to 0.05. We find that $T_{n}$ performs equally well in dimensional control in various situations and can be controlled at a empirical sizes of about 0.05. \\
The empirical powers of our proposed test statistic and the other three comparison statistics are displayed in Figures 7-12. Again we only report the empirical powers of $T_{n}$ under N(0,1), the remaining two distributions are shown in the Supplementary Material. Its empirical powers are characterized similarly to Simulation \ref{sec3.1}, with powers stronger than the other three tests at $\rho=0.1$ or $\rho=0.2$ and lower for all tests at $\rho=0.3$ or $\rho=0.4$.

\begin{table}[H]
\begin{center}
\caption{{Empirical sizes in the GLHT problem when $z_{ij}$ follows $N(0,1)$ for scenario 1.}}
\label{tab3}
\begin{tabularx}{\textwidth}{XXXXXX}
\hline
$p$  &$n$  &\multicolumn{4}{c}{$\mathbf{z}_{ijl}\overset{i.i.d}{\sim}N(0,1)$} \\ 
\cmidrule(r){3-6} 
                    &  &$T_{n}$ &$T_{zb}$  &$T_{zjt}$ &$T_{zz}$  \\
\hline
100    &$n_{1}$  &0.0516	&0.0532	&0.0478	&0.0458	  \\
       &$n_{2}$  &0.0504	&0.0566	&0.0462	&0.0448	    \\
       &$n_{3}$  &0.0528	&0.062	  &0.0558	&0.0546	  \\
200    &$n_{1}$  &0.0516	&0.053	  &0.0522	&0.0504	  \\
       &$n_{2}$  &0.0502	&0.0614	&0.0528	&0.051	 \\
       &$n_{3}$  &0.0496	&0.0582	&0.0512	&0.0496	   \\
300    &$n_{1}$  &0.0508	&0.0532	&0.05	  &0.0488	    \\
       &$n_{2}$  &0.0514	&0.0558	&0.049	  &0.047	  \\
       &$n_{3}$  &0.0494	&0.052	  &0.0466	&0.0454   \\
400    &$n_{1}$  &0.0532	&0.0476	&0.0438	&0.043	   \\
       &$n_{2}$  &0.0514	&0.0542	&0.0464	&0.0454	 \\
       &$n_{3}$  &0.0554	&0.066	  &0.0608	&0.0598	  \\
500    &$n_{1}$  &0.0518	&0.0532	&0.0496	&0.0482	      \\
       &$n_{2}$  &0.0526	&0.0598	&0.0578	&0.0568	  \\
       &$n_{3}$  &0.0494	&0.053	  &0.0488	&0.0478	\\
\hline                                                                                                                     
\end{tabularx}
\end{center}
\end{table}

 \begin{figure}[H]
  \centering
  \hspace{-0.8cm}
  \includegraphics[scale=1]{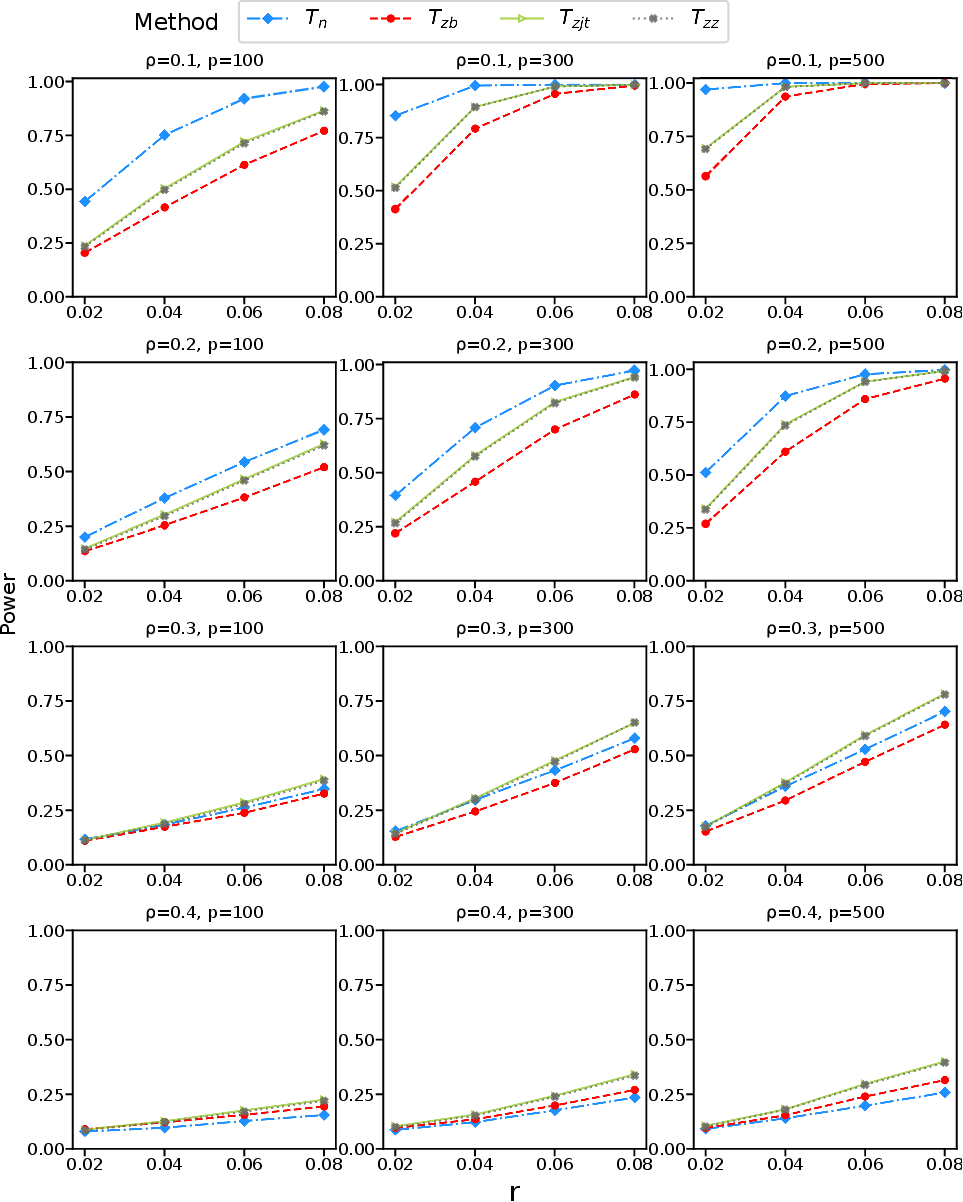}
  \begin{flushleft}
    {\par{\footnotesize\textbf{Figure 7.} Empirical powers when $z_{ij}$ follows $N(0,1)$ and $n_{1}=(25,30,40,45)$ for scenario 1 under different signal levels of $r$ and sparsity levels of $\rho$ in GLHT problem.}}
    \end{flushleft}
 \end{figure}

 \begin{figure}[H]
  \centering
  \hspace{-0.8cm}
   \includegraphics[scale=1]{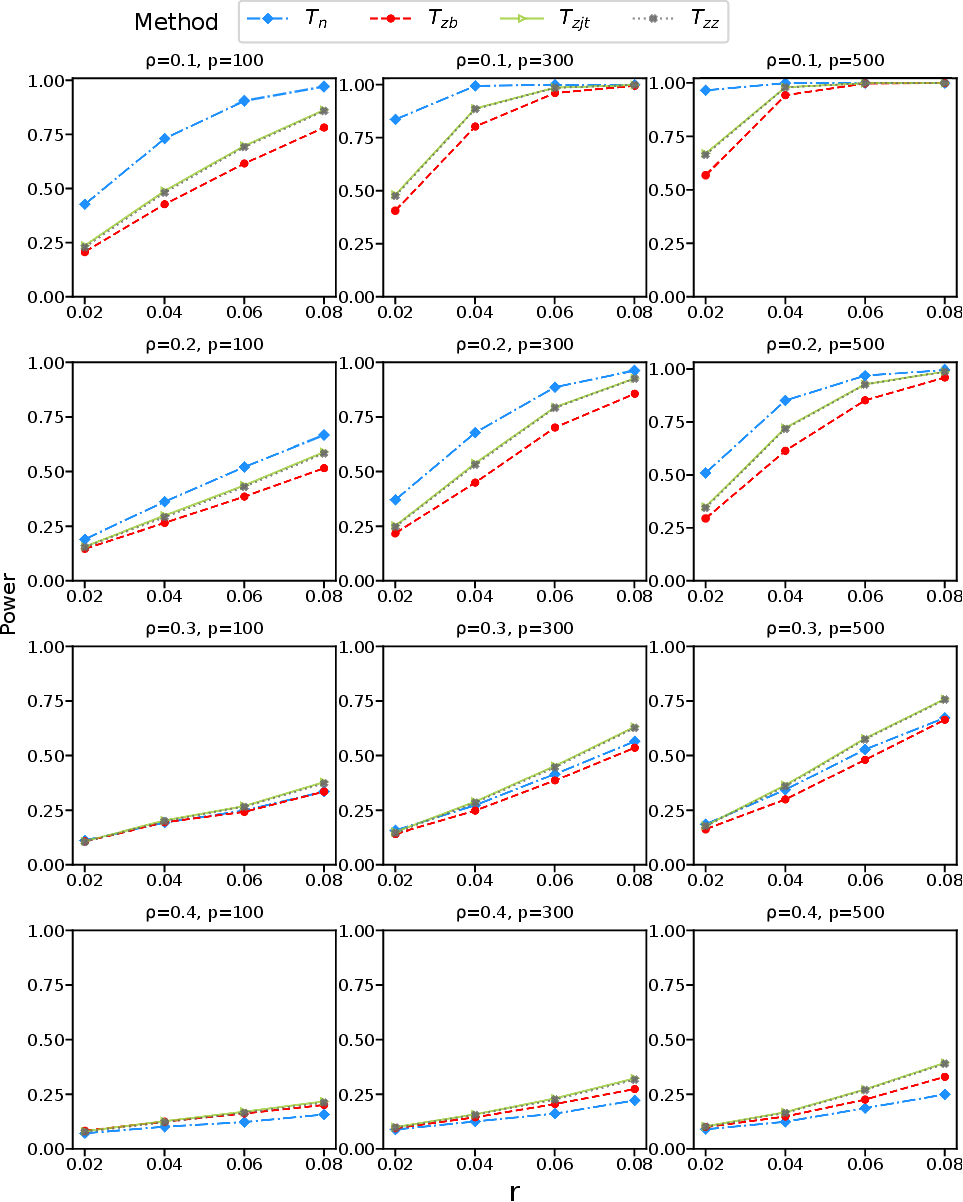}
   \begin{flushleft}
    {\par{\footnotesize\textbf{Figure 8.} Empirical powers when $z_{ij}$ follows $N(0,1)$ and $n_{2}=(30,40,55,60)$ for scenario 1 under different signal levels of $r$ and sparsity levels of $\rho$ in GLHT problem.}}
    \end{flushleft}
 \end{figure}
 
  \begin{figure}[H]
   \centering
   \hspace{-0.8cm}
   \includegraphics[scale=1]{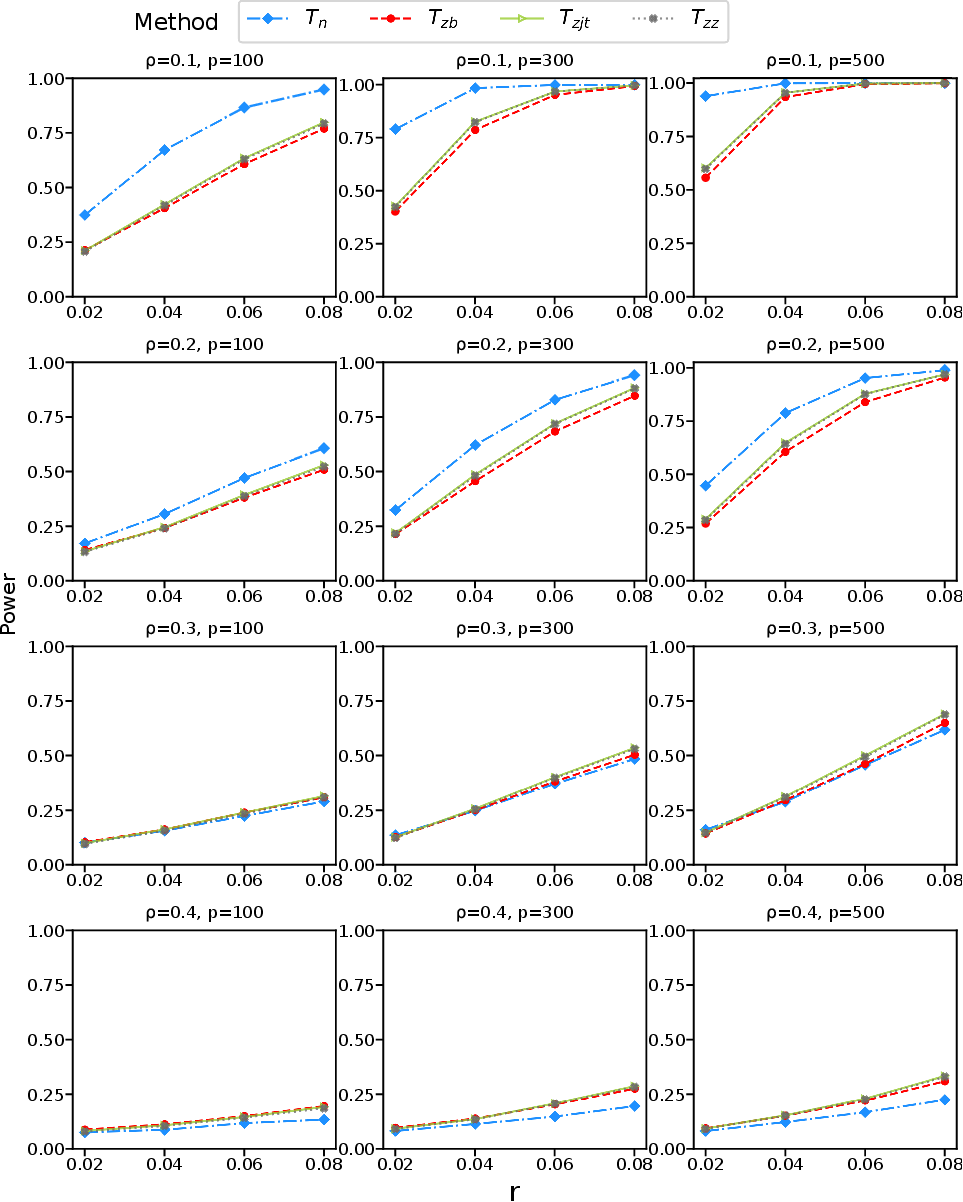}
   \begin{flushleft}
    {\par{\footnotesize\textbf{Figure 9.} Empirical powers when $z_{ij}$ follows $N(0,1)$ and $n_{3}=(80,95,115,120)$ for scenario 1 under different signal levels of $r$ and sparsity levels of $\rho$ in GLHT problem.}}
    \end{flushleft}
 \end{figure}
 
\begin{table}[H]
\begin{center}
\caption{{Empirical sizes in the GLHT problem when $z_{ij}$ follows $N(0,1)$ for scenario 2.}}
\label{tab4}
\begin{tabularx}{\textwidth}{XXXXXX}
\hline
$p$  &$n$  &\multicolumn{4}{c}{$\mathbf{z}_{ijl}\overset{i.i.d}{\sim}N(0,1)$} \\ 
\cmidrule(r){3-6} 
                    &  &$T_{n}$ &$T_{zb}$  &$T_{zjt}$ &$T_{zz}$  \\
\hline
100    &$n_{1}$   &0.0504	&0.056	  &0.0542	&0.0534 \\    
       &$n_{2}$   &0.0512	&0.0574	&0.052	  &0.0518	    \\    
       &$n_{3}$   &0.0512	&0.0544	&0.048	  &0.048	  \\    
200    &$n_{1}$   &0.05	  &0.0528	&0.05	  &0.05	      \\    
       &$n_{2}$   &0.0528	&0.0586	&0.0536	&0.0534	     \\    
       &$n_{3}$   &0.0512	&0.0556	&0.0514	&0.0514	   \\    
300    &$n_{1}$   &0.0482	&0.0526	&0.0488	&0.0488	   \\    
       &$n_{2}$   &0.0504	&0.0554	&0.0502	&0.049	     \\    
       &$n_{3}$   &0.0494	&0.048	  &0.0456	&0.0454  \\    
400    &$n_{1}$   &0.0498	&0.0536	&0.0508	&0.0498     \\    
       &$n_{2}$   &0.0502	&0.0544	&0.048	  &0.047	    \\    
       &$n_{3}$   &0.0498	&0.0536	&0.0508	&0.0498	      \\    
500    &$n_{1}$   &0.0516	&0.051	  &0.0484	&0.0478	  \\    
       &$n_{2}$   &0.0502	&0.0504	&0.046	  &0.0456	   \\    
       &$n_{3}$   &0.0546	&0.0586	&0.0566	&0.0564	   \\    
\hline                                                                                                                           
\end{tabularx}
\end{center}
\end{table}
 
  \begin{figure}[H]
   \centering
   \hspace{-0.8cm}
   \includegraphics[scale=1]{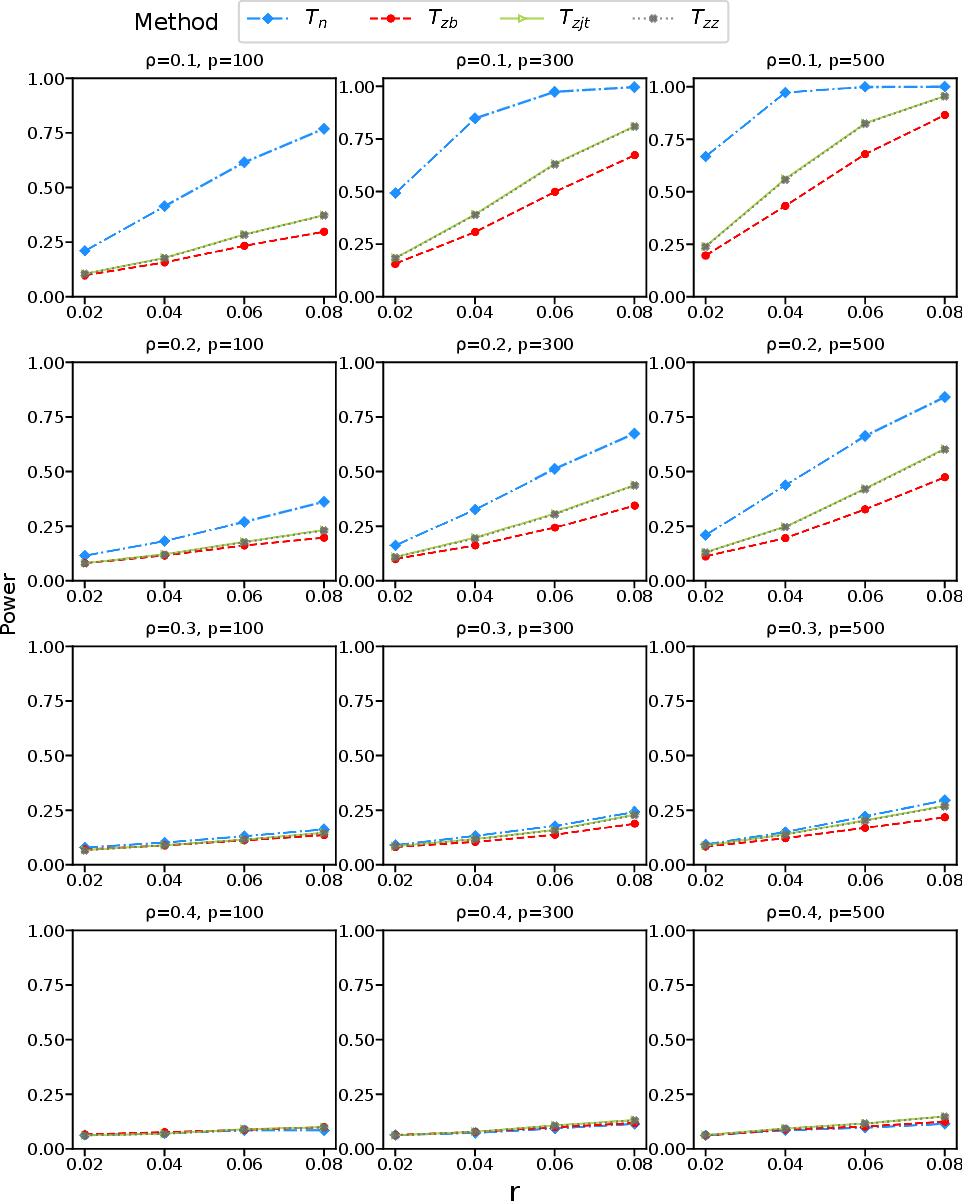}
   \begin{flushleft}
    {\par{\footnotesize\textbf{Figure 10.} Empirical powers when $z_{ij}$ follows $N(0,1)$ and $n_{1}=(25,30,40,45)$ for scenario 2 under different signal levels of $r$ and sparsity levels of $\rho$ in GLHT problem.}}
    \end{flushleft}
 \end{figure}
 
  \begin{figure}[H]
   \centering
   \hspace{-0.8cm}
   \includegraphics[scale=1]{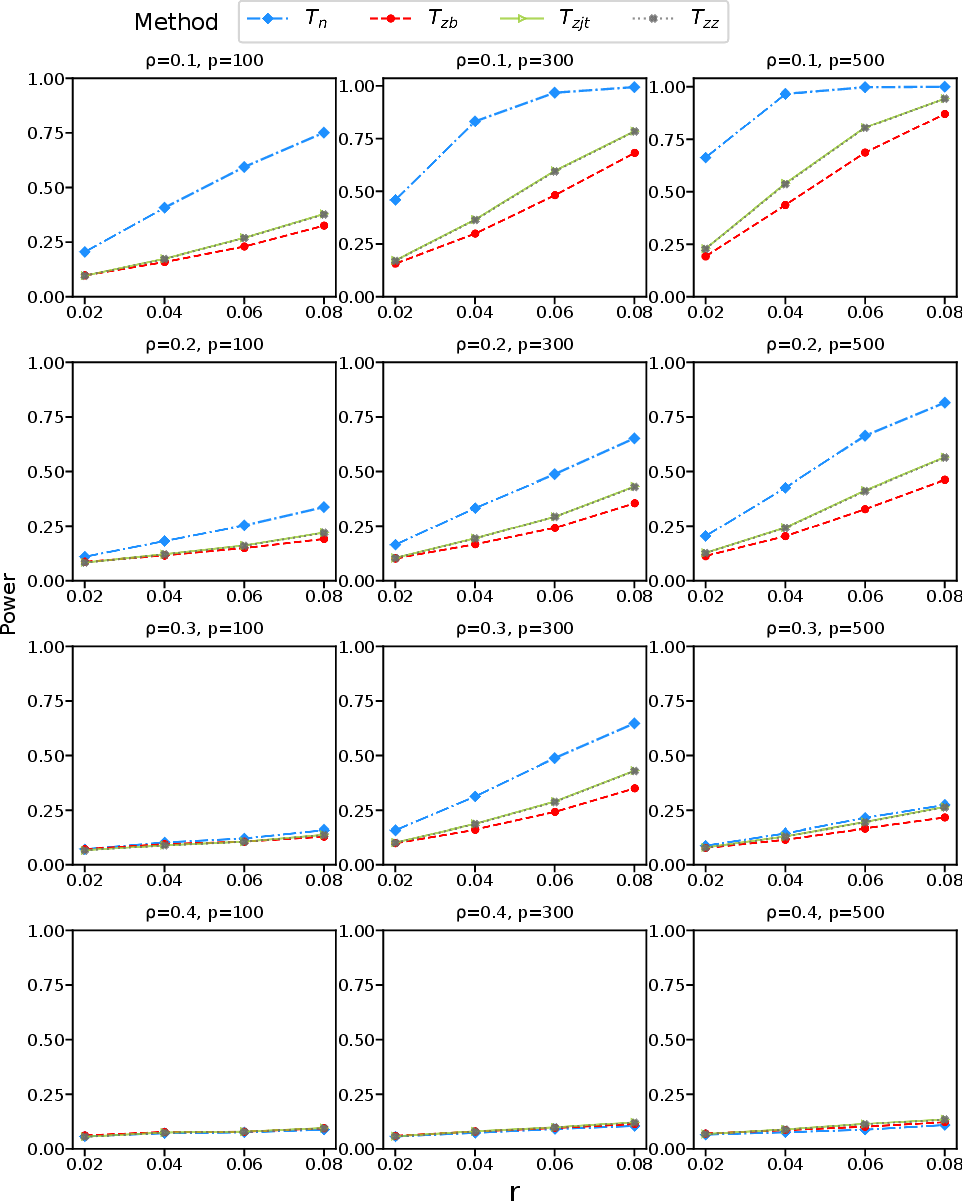}
   \begin{flushleft}
    {\par{\footnotesize\textbf{Figure 11.} Empirical powers when $z_{ij}$ follows $N(0,1)$ and $n_{2}=(30,40,55,60)$ for scenario 2 under different signal levels of $r$ and sparsity levels of $\rho$ in GLHT problem.}}
    \end{flushleft}
 \end{figure}
 
  \begin{figure}[H]
   \centering
   \hspace{-0.8cm}
    \includegraphics[scale=1]{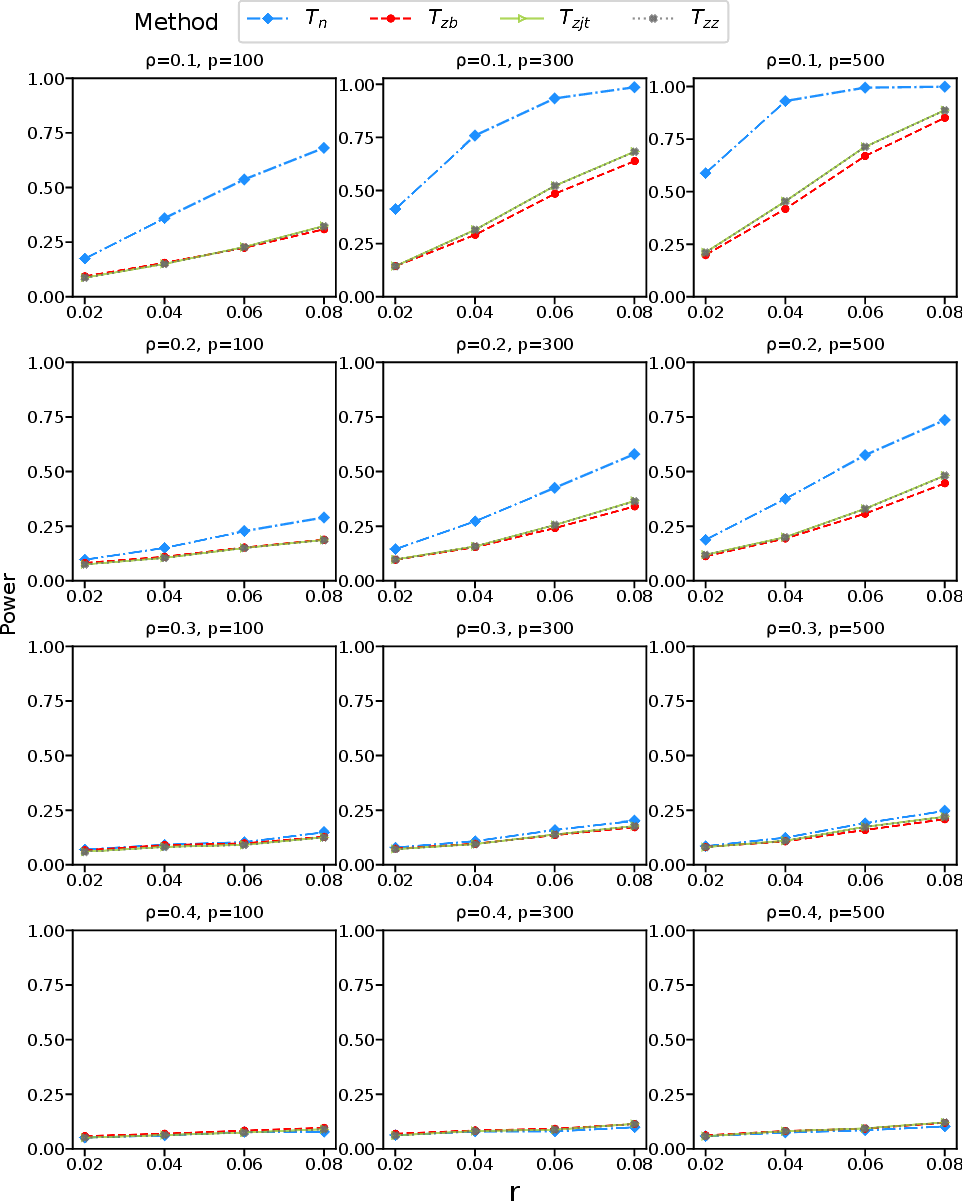}
    \begin{flushleft}
    {\par{\footnotesize\textbf{Figure 12.} Empirical powers of when $z_{ij}$ follows $N(0,1)$ and $n_{3}=(80,95,115,120)$ for scenario 2 under different signal levels of $r$ and sparsity levels of $\rho$ in GLHT problem.}}
    \end{flushleft}
 \end{figure}
\section{Real data analysis}\label{sec4}
\noindent 
In this section, we analyze a real data example, and the data used is corneal data, which comes from Locantore et al.(1999). It is divided into four groups of corneal surface data, the first group consists of 43 healthy corneal surface data and is called the normal cornea group, while the other three groups 
are the abnormal cornea groups, which are the unilateral suspected, suspicious map, and clinical keratoconus cornea groups, which are each composed of 14, 21, and 72 corneal surface data reflecting varying degrees of keratoconus cornea, which is a condition that occurs when the cornea is misshapen. We denote these four sets of data as 'Nor', 'Uni', 'Sus' and 'Ker', respectively. The entire reconstructed corneal surface data set is represented by 150 vectors of 
length 2000.\\
Without assuming that the four sets of corneal data have the same covariance array, the test we are interested in is whether the four sets of corneal data 
have the same mean. Since the total sample size of the data set is n=150 and the dimension p=2000, the statistic $T_{n}$ proposed in this paper and the statistics $T_{zb}$, $T_{zjt}$ and $T_{hj}$ proposed in \cite{Zhou et al:2017}, \cite{Zhang et al:2022(b)} and \cite{Hu et al:2017} are applied, respectively. The results of the test for the mean of the corneal surface data for the four groups are in the first row of Table \ref{tab5}. It can be seen that all 
p-values are very small, so all tests indicate a very strong rejection of the null hypothesis, suggesting that the means of the four sets of corneal data are unlikely to be the same. Further, it is also seen that the degrees of freedom of approximation for the estimates of $T_{n}$ and $T_{zjt}$ are not large, suggesting that the corneal data may be moderately or highly correlated, so that the normal approximation used for $T_{zb}$ and $T_{hj}$ are not sufficient to approximate a correlated null distribution, even though the p-values for all tests are very small. To test whether the highly significant results were affected by sample size, we next examined whether any two corneal groups had different mean corneal surface data. We again apply the four statistics to the 
two-sample problem, and in Table \ref{tab5} we can see that the four statistics give similar conclusions for the two-sample tests listed in the table. All 
the test results show that our proposed test statistic is extremely significant on any pair of corneal data mean detection, indicating a strong rejection of the null hypothesis that the four corneal data means are not equal.

\begin{table}[H]
\begin{center} 
\scriptsize   
\caption{{Results of $T_{n}, T_{zb}, T_{zjt}$ and $T_{hj}$ for testing the mean surface differences between the four groups corneal data.}} 
\label{tab5}                                                                                  
\begin{tabularx}{\textwidth}{XXXXXX}                                                  
\hline                                                                                                                                                 
&Method        &Statistic  &P-value &$\hat{\beta}$  &$\hat{d}$ \\  
\hline                                                           
All Groups    &$T_{n}$   &2698.2660  &0.0000  &41.1986  &7.5304  \\                         
              &$T_{zb}$      &202.2157  &0.0000  &-  &-  \\                      
              &$T_{zjt}$          &159.7330 &0.0000   &6.1464 &6.1652 \\   
              &$T_{hj}$          &3.0465 &0.0012   &- &- \\               
Uni vs. Sus   &$T_{n}$   &123.1895  &0.0000  &40.8643  &2.6409\\                         
              &$T_{zb}$      &20.3225  &0.0000  &-  &- \\                      
              &$T_{zjt}$          &5.9080 &$1.7313\times10^{-9}$   &6.7825    &2.0451\\  
              &$T_{hj}$          &-0.4375 &0.6691   &- &- \\                
Uni vs. Nor   &$T_{n}$   &33.1273  &0.0000  &45.5852  &2.7707 \\                         
              &$T_{zb}$      &109.2432  &0.0000  &-  &-  \\                      
              &$T_{zjt}$          &5.1537 &$1.2767\times10^{-7}$   &7.5138  &2.1161 \\       
              &$T_{hj}$          &-0.5106 &0.6952   &- &- \\           
Uni vs. Ker   &$T_{n}$   &572.1285  &0.0000  &52.7285  &2.7070  \\                         
              &$T_{zb}$      &180.7023  &0.0000  &-  &- \\                      
              &$T_{zjt}$          &32.2614 &0.0000   &8.1939  &2.0923  \\   
              &$T_{hj}$          &0.5692 &0.2846   &- &- \\               
Sus vs. Nor   &$T_{n}$   &283.2604  &0.0000  &23.6926  &2.6193  \\                         
              &$T_{zb}$      &95.9070  &0.0000  &-  &-  \\                      
              &$T_{zjt}$          &19.2810 &0.0000   &3.9590  &2.1440  \\   
              &$T_{hj}$          &0.7712 &0.2203   &- &- \\
Sus vs. Ker   &$T_{n}$   &1868.0213  &0.0000  &26.2497  &2.8569  \\           
               &$T_{zb}$      &246.7162  &0.0000  &-  &-  \\                   
              &$T_{zjt}$          &95.9644 &0.0000   &3.9744  &2.2941  \\    
              &$T_{hj}$          &5.8474 &$2.4959\times10^{-9}$   &- &- \\            
Nor vs. Ker   &$T_{n}$   &1369.4114  &0.0000  &33.2393  &3.0658  \\            
               &$T_{zb}$      &83.3126  &0.0000  &-  &-  \\                   
              &$T_{zjt}$          &91.5327 &0.0000   &4.2381  &2.7170  \\                
              &$T_{hj}$          &4.9935 &$2.9642\times10^{-7}$   &- &- \\               
\hline                                                                                                                                                
\end{tabularx}                                                                                                                                        
\end{center}                                                                                                                                          
\end{table}

\section{Concluding remarks}\label{sec5}
\noindent Various methods have been developed to deal with the linear hypothesis testing problem in high-dimensional environments. Currently, this problem is still an active research topic in statistics. Although meaningful progress has been made, further research is still needed. In this paper, we propose a random integration method to study the linear hypothesis testing problem by modifying the test statistic proposed in \cite{Zhang et al:2022(b)} based on the $L^{2}$-norm, and the new test statistic performs better in some cases. For example, when nonzero signals are weakly dense with nearly the same sign or when there are more dense or only weakly dense nonzero signals, our proposed test performs well in handling such problems. For the GLHT problem, the proposed test statistic and the chi-square mixture have the same normal or non-normal limiting distribution under certain conditions, it is then suggested that the distribution of the chi-square type mixture be used to approximate the null distribution of the test. We establish the asymptotic power of this test under local alternative conditions and study the effect of data non-normality. Several simulation studies demonstrate good performance in terms of dimensional control and the power of our test. \\
It is worth noting that there are further issues to be investigated with our proposed method. Firstly, in our simulation studies and real data analysis, we use p-dimensionally independent density functionals with the limited $\alpha_{i}$ and $\beta_{i}$ as the weight functions. Therefore, an interesting future topic is to consider the choice of a p-dimensional independent density function or other weighting function. Secondly, the test remains reliable when there are many small to moderate component differences or when non-zero signals have nearly identical signals. We are given a weight function and a density function with independent components. In further work, we will consider a different metric with dependent components. Finally, we will consider the other settings in more detail to ensure that they meet the requirements of the actual application.

\section*{Acknowledgments}

Dr Cao's research is supported by Humanities and Social Sciences Fund of the Ministry of Education (No. 22YJC910001).

\section*{Supplementary Materials}
\noindent The results of the simulations under the two remaining distributions are shown in the Supplementary Material.
\titleformat{\section}[hang]{\normalfont\Large\bfseries}{Appendix}{0pt}{}
\newcommand{\proof}{{{\it Proof.}}}
\begin{appendices}
\section{}
\textbf{Proof of Theorem \ref{th2.1}}\\
\noindent The test statistic $T_{n}$ obtained by the random integration method can be expressed in the following form,
\begin{equation*}
\begin{aligned}
T_{n}&=\|(I_{q}\otimes\delta^{T})C\hat{\mu}\|^{2}\\
&=\int(C\hat{\mu})^{T}(I_{q}\otimes\delta\delta^{T})(C\hat{\mu})w(\delta)d\delta\\
&=\int\left(\sum_{i=1}^{qp}c_{i}\hat{\mu}_{i}\delta_{i}\right)^{2}w(\delta)d\delta\\
&=\int\left(\sum_{i=1}^{qp}c_{i}^{2}\hat{\mu}_{i}^{2}\delta_{i}^{2}+\sum_{i{\neq}j}^{qp}c_{i}c_{j}\hat{\mu}_{i}\hat{\mu}_{j}\delta_{i}\delta_{j}\right)w(\delta)d\delta\\
&=\sum_{i=1}^{qp}c_{i}^{2}\hat{\mu}_{i}^{2}\int\delta_{i}^{2}w(\delta)d\delta+\sum_{i{\neq}j}^{qp}c_{i}c_{j}\hat{\mu}_{i}\hat{\mu}_{j}\int\delta_{i}\delta_{j}w(\delta)d\delta\\
&=\sum_{i=1}^{qp}c_{i}^{2}\hat{\mu}_{i}^{2}\sum_{i=1}^{p}q(\beta_{i}^{2}+\alpha_{i}^{2})+\sum_{i{\neq}j}^{qp}c_{i}c_{j}\hat{\mu}_{i}\hat{\mu}_{j}\sum_{i{\neq}j}^{p}q\alpha_{i}\alpha_{j}\\
&=(C\hat{\mu})^{T}(I_{q}{\otimes}W_{p})(Cu).
\end{aligned}
\end{equation*}
where $c_{i}$ and $\mu_{i}$ are elements of $C$ and $\mu$, respectively. $w(\delta_{i})$ is a density function with the mean $\alpha_{i}$ and the variance $\beta_{i}^{2}$ for $i=(1,\ldots,p)$, and denotes $W_{p}=B+aa^{T}$, where $a=(\alpha_{1},\ldots,\alpha_{p})^{T}$,

$\\$
$B=\begin{pmatrix}
\beta_{1}^{2} & 0 & \cdots & 0 \\
0 & \beta_{2}^{2} & \cdots & 0 \\
\vdots & \vdots & \ddots & \vdots \\
0 & 0 & \cdots & \beta_{p}^{2}
\end{pmatrix}$.\\
~\\
Then, we have
\begin{equation*}
\begin{aligned}
T_{n}&=(C\hat{\mu})^{T}(I_{q}{\otimes}W_{p})(Cu)\\
&=\hat{\mu}^{T}(G^{T}{\otimes}I_{p})(I_{q}{\otimes}W_{p})(G{\otimes}I_{p})\hat{\mu}\\
&=\hat{\mu}^{T}(G^{T}G{\otimes}W_{p})\hat{\mu}\\
&=\hat{\mu}^{T}(H{\otimes}W_{p})\hat{\mu}.
\end{aligned}
\end{equation*}
This completes the proof of Theorem \ref{th2.1}.\\

\textbf{Proof of Theorem \ref{th2.2}}\\
\noindent After the transformation, we get $T_{n,0}=\varphi^{T}(A{\otimes}W_{p})\varphi$, when all k samples are normal, we get $\varphi{\sim}N(0,\Sigma)$. We denote $Z=\Sigma^{-1/2}\varphi{\sim}N(0,I_{kp})$, we have
\begin{equation*}
\begin{aligned}
T_{n,0}&=Z^{T}\Sigma^{1/2}(A{\otimes}W_{p})\Sigma^{1/2}Z\\
&=Z^{T}Q^{T}{\Lambda}QZ\\
&=U^{T}{\Lambda}U\\
&=\sum_{r=1}^{kp}\lambda_{n,r}A_{r},
\end{aligned}
\end{equation*}
where $\lambda_{n,r}$ and $A_{r}$ are elements of $\Lambda$ and $U^{T}U$, we make an eigenvalue transformation of $\Sigma^{1/2}(A{\otimes}W_{p})\Sigma^{1/2}$,
\begin{equation*}
\begin{aligned}
\Sigma^{1/2}(A{\otimes}W_{p})\Sigma^{1/2}&\rightarrow\Sigma(A{\otimes}W_{p})\\
&\rightarrow\Sigma(B^{T}{\otimes}W_{P})(B{\otimes}W_{P})\\
&\rightarrow(B{\otimes}W_{P})\Sigma(B^{T}{\otimes}W_{P}),
\end{aligned}
\end{equation*}
then, we have $\Sigma^{1/2}(A{\otimes}W_{p})\Sigma^{1/2}=Q^{T}{\Lambda}Q$, and $U=QZ{\sim}N(0,1)$, where $\Lambda$ and $Q$ are the eigenvalues and eigenvectors of $\Sigma^{1/2}(A{\otimes}W_{p})\Sigma^{1/2}$, respectively. Then, we denote $T_{n,0}^{*}=\sum_{r=1}^{kp}\lambda_{n,r}A_{r}$, this completes the proof of Theorem \ref{th2.2}.\\

\textbf{Proof of Theorem \ref{th2.3}}\\
\noindent We prove the first expression of \ref{eq2.9} by the characteristic function($\psi_{X}(t)=E(e^{itX}$) for a random variable $X$) method. Denotes $\varphi_{n,p}=(B{\otimes}W_{p})\varphi$ where $B$ and $\varphi$ is given in the preceding text. By some algebraic calculations, we have $E(\varphi_{n,p})=0$ and $Cov(\varphi_{n,p})=\Omega_{n}$, where $\Omega_{n}$ is a $(kp{\times}kp)$ semi-positive matrix. Denotes $u_{n,r}$ are the eigenvectors associated with the decreasing-ordered eigenvalues $\lambda_{n,r}$ of $\Omega_{n}$, where $r\in\{1,\ldots,kp\}$. Then, we have $\varphi_{n,p}=\sum_{r=1}^{kp}\eta_{n,r}u_{n,r}$, where $\eta_{n,r}=\varphi_{n,p}^{T}u_{n,r}$, $r\in\{1,\ldots,kp\}$, and we get $E(\eta_{n,r})=0$ and $Var(\eta_{n,r})=\lambda_{n,r}$, $r\in\{1,\ldots,kp\}$. Denotes $v_{r}=(B{\otimes}W_{p})u_{n,r}=(v_{1,r}^{T},\ldots,v_{k,r}^{T})^{T}$, $r\in\{1,\ldots,kp\}$, and $v_{i,r},~i\in\{1,\ldots,k\}$ are $p{\times}1$ vectors. Then we have $\eta_{n,r}=\varphi_{n,p}^{T}u_{n,r}=\varphi^{T}(B^{T}{\otimes}W_{p})u_{n,r}=\sum_{i=1}^{k}v_{i,r}^{T}\varphi_{i}$. Thus, we have
\begin{equation*}
\begin{aligned}
\eta_{n,r}^{2}&=\left(\sum_{i=1}^{k}v_{i,r}^{T}\varphi_{i}\right)^{2}\\
&=\sum_{i=1}^{k}(v_{i,r}^{T}\varphi_{i})^{2}+2\sum_{1{\leq}i<j{\leq}k}(v_{i,r}^{T}\varphi_{i})v_{j,r}^{T}\varphi_{i}.
\end{aligned}
\end{equation*}
We know that $E(\varphi_{i})=0$ and $Cov(\varphi_{i})=\Sigma_{i}$, $i\in\{1,\ldots,k\}$ and $\varphi_{i}$ are independent, we have
\begin{equation*}
\begin{aligned}
Var(\eta_{n,r}^{2})&=\sum_{i=1}^{k}Var\{(v_{i,r}^{T}\varphi_{i})^{2}\}+4\sum_{1{\leq}i<j{\leq}k}Var\{(v_{i,r}^{T}\varphi_{i})(v_{j,r}^{T}\varphi_{i})\},
\end{aligned}
\end{equation*}
The following calculations are performed separately for $Var\{(v_{i,r}^{T}\varphi_{i})^{2}\}$ and $Var\{(v_{i,r}^{T}\varphi_{i})(v_{j,r}^{T}\varphi_{i})\}$, notice that $\varphi_{i}=\sqrt{n_{i}}(\bar{y_{i}}-\mu_{i})$, denotes $x_{i}=y_{i}-\mu_{i}$, thus, we have $E(x_{i})=0,~Cov(x_{i})=\Sigma_{i}$, and $\bar{y_{i}}-\mu_{i}=\frac{1}{n_{i}}\sum_{i=1}^{n_{i}}y_{i}-\mu_{i}=\frac{1}{n_{i}}\sum_{i=1}^{n_{i}}x_{i}$. Then, we have
\begin{equation*}
\begin{aligned}
Var\{(v_{i,r}^{T}\varphi_{i})^{2}\}&=Var\{(v_{i,r}^{T}\sqrt{n_{i}}(\bar{y_{i}}-\mu_{i}))^{2}\}\\
&=Var\{(v_{i,r}^{T}\frac{\sqrt{n_{i}}}{n_{i}}\sum_{i=1}^{n_{i}}x_{i})^{2}\}\\
&=n_{i}^{-2}Var\{(v_{i,r}^{T}\sum_{i=1}^{n_{i}}x_{i})^{2}\}\\
&=n_{i}^{-2}Var\{v_{i,r}^{T}(\sum_{i=1}^{n_{i}}||x_{i}||^{2})v_{i,r}+2v_{i,r}^{T}(\sum_{1{\leq}i<j{\leq}n_{i}}x_{i}^{T}x_{j})v_{i,r}\}\\
&=n_{i}^{-2}\{(v_{i,r}^{T}v_{i,r})^{2}\sum_{i=1}^{n_{i}}Var(||x_{i}||^{2})+2(v_{i,r}^{T}v_{i,r})^{2}\sum_{i{\neq}j}Var(x_{i}^{T}x_{j})\}\\
&=n_{i}^{-2}\{(v_{i,r}^{T}v_{i,r})^{2}n_{i}Var(||x_{1}||^{2})+2(v_{i,r}^{T}v_{i,r})^{2}n_{i}(n_{i}-1)\Sigma_{i}^{2}\}\\
&=n_{i}^{-2}\{(v_{i,r}^{T}v_{i,r})^{2}n_{i}[E(||x_{1}||^{4})-E^{2}(||x_{1}||^{2})]+2(v_{i,r}^{T}v_{i,r})^{2}n_{i}(n_{i}-1)\Sigma_{i}^{2}\}\\
&=n_{i}^{-2}\{(v_{i,r}^{T}v_{i,r})^{2}n_{i}[E(y_{i1}-\mu_{i})^{4}-\Sigma_{i}^{2}]+2(v_{i,r}^{T}v_{i,r})^{2}n_{i}(n_{i}-1)\Sigma_{i}^{2}\}\\
&=n_{i}^{-1}\{E[\{v_{i,r}^{T}(y_{i1}-\mu_{i})\}^{4}]-(v_{i,r}^{T}\Sigma_{i}v_{i,r})^{2}\}+2(v_{i,r}^{T}\Sigma_{i}v_{i,r})^{2}-2n_{i}^{-1}(v_{i,r}^{T}\Sigma_{i}v_{i,r})^{2}\\
&=2(v_{i,r}^{T}\Sigma_{i}v_{i,r})^{2}+\frac{E[\{v_{i,r}^{T}(y_{i1}-\mu_{i})\}^{4}]-3(v_{i,r}^{T}\Sigma_{i}v_{i,r})^{2}}{n_{i}}.\\
\end{aligned}
\end{equation*}
Similarly,
\begin{equation*}
\begin{aligned}
Var\{(v_{i,r}^{T}\varphi_{i})(v_{j,r}^{T}\varphi_{i})\}&=Var\{(v_{i,r}^{T}\frac{\sqrt{n_{i}}}{n_{i}}\sum_{i=1}^{n_{i}}x_{i})(v_{j,r}^{T}\frac{\sqrt{n_{j}}}{n_{j}}\sum_{j=1}^{n_{j}}x_{j})\}\\
&=\frac{1}{n_{i}n_{j}}Var\{(v_{i,r}^{T}\sum_{i=1}^{n_{i}}x_{i})(v_{j,r}^{T}\sum_{j=1}^{n_{j}}x_{j})\}\\
&=\frac{1}{n_{i}n_{j}}\{E\{[(v_{i,r}^{T}\sum_{i=1}^{n_{i}}x_{i})(v_{j,r}^{T}\sum_{j=1}^{n_{j}}x_{j})]^{2}\}-E^{2}\{(v_{i,r}^{T}\sum_{i=1}^{n_{i}}x_{i})(v_{j,r}^{T}\sum_{j=1}^{n_{j}}x_{j})\}\}\\
&=\frac{1}{n_{i}n_{j}}E\{(v_{i,r}^{T}\sum_{i=1}^{n_{i}}x_{i})^{2}\}E\{(v_{j,r}^{T}\sum_{j=1}^{n_{j}}x_{j})^{2}\}\\
&=\frac{1}{n_{i}n_{j}}n_{i}v_{i,r}^{T}\Sigma_{i}v_{i,r}n_{j}v_{j,r}^{T}\Sigma_{j}v_{j,r}\\
&=v_{i,r}^{T}\Sigma_{i}v_{i,r}v_{j,r}^{T}\Sigma_{j}v_{j,r}.
\end{aligned}
\end{equation*}
By applying proposition A.1 (i) in \cite{Chen et al:2010+}, we have the following result
\begin{equation*}
\begin{aligned}
E[\{v_{i,r}^{T}(y_{i1}-\mu_{i})\}^{4}]&=E\{(v_{i,r}^{T}\Gamma_{i}z_{ij})^{4}\}\\
&=E\{(z_{ij}^{T}\Gamma_{i}^{T}v_{i,r}v_{i,r}^{T}\Gamma_{i}z_{ij})^{2}\}\\
&=tr^{2}(\Gamma_{i}^{T}v_{i,r}v_{i,r}^{T}\Gamma_{i})+2tr\{(\Gamma_{i}^{T}v_{i,r}v_{i,r}^{T}\Gamma_{i})^{2}\}+{\Delta}tr(\Gamma_{i}^{T}v_{i,r}v_{i,r}^{T}\Gamma_{i}\circ\Gamma_{i}^{T}v_{i,r}v_{i,r}^{T}\Gamma_{i})\\
&=3(\Gamma_{i}^{T}v_{i,r}v_{i,r}^{T}\Gamma_{i})^{2}+{\Delta}tr(\Gamma_{i}^{T}v_{i,r}v_{i,r}^{T}\Gamma_{i}\circ\Gamma_{i}^{T}v_{i,r}v_{i,r}^{T}\Gamma_{i})\\
&~{\leq}(3+\Delta)(\Gamma_{i}^{T}v_{i,r}v_{i,r}^{T}\Gamma_{i})^{2}\\
&=(3+\Delta)tr(v_{i,r}^{T}\Gamma_{i}\Gamma_{i}^{T}v_{i,r})^{2}\\
&=(3+\Delta)(v_{i,r}^{T}\Sigma_{i}v_{i,r})^{2},
\end{aligned}
\end{equation*}
where $\circ$ denoting the Hardmard product operator. Thus, we have
\begin{equation*}
\begin{aligned}
Var(\eta_{n,r}^{2})&=\sum_{i=1}^{k}Var\{(v_{i,r}^{T}\varphi_{i})^{2}\}+4\sum_{1{\leq}i<j{\leq}k}Var\{(v_{i,r}^{T}\varphi_{i})(v_{j,r}^{T}\varphi_{i})\}\\
&~{\leq}\sum_{i=1}^{k}(2+\Delta/n_{i})(v_{i,r}^{T}\Sigma_{i}v_{i,r})^{2}+4\sum_{1{\leq}i<j{\leq}k}v_{i,r}^{T}\Sigma_{i}v_{i,r}v_{j,r}^{T}\Sigma_{j}v_{j,r}\\
&=2\sum_{i=1}^{k}\sum_{j=1}^{k}v_{i,r}^{T}\Sigma_{i}v_{i,r}v_{j,r}^{T}\Sigma_{j}v_{j,r}+\sum_{i=1}^{k}(v_{i,r}^{T}\Sigma_{i}v_{i,r})^{2}\Delta/n_{i}\\
&~{\leq}2(\sum_{i=1}^{k}v_{i,r}^{T}\Sigma_{i}v_{i,r})^{2}+\Delta/n_{min}\sum_{i=1}^{k}(v_{i,r}^{T}\Sigma_{i}v_{i,r})^{2}\\
&~{\leq}(2+\Delta/n_{min})(\sum_{i=1}^{k}v_{i,r}^{T}\Sigma_{i}v_{i,r})^{2}\\
&=(2+\Delta/n_{min})\lambda_{n,r},
\end{aligned}
\end{equation*}
where $n_{min}=min_{i=1}^{n}n_{i}$ and we have
\begin{equation*}
\begin{aligned}
\sum_{i=1}^{k}v_{i,r}^{T}\Sigma_{i}v_{i,r}=v_{r}^{T}{\Sigma}v_{r}=u_{n,r}^{T}(B{\otimes}W_{p}){\Sigma}(B^{T}{\otimes}W_{p})u_{n,r}=\lambda_{n,r}.
\end{aligned}
\end{equation*}
Denotes $T_{n,0}=\sum_{r=1}^{kp}\eta_{n,r}^{2}$, we set
\begin{equation*}
\begin{aligned}
&\tilde{T}_{n,0}=\{T_{n,0}-tr(\Omega_{n})\}/\sqrt{2tr(\Omega_{n}^{2})}=\sum_{r=1}^{kp}(\eta_{n,r}^{2}-\lambda_{n,r})/\sqrt{2tr(\Omega_{n}^{2})}\\
&\tilde{T}_{n,0}^{q}=\sum_{r=1}^{q}(\eta_{n,r}^{2}-\lambda_{n,r})/\sqrt{2tr(\Omega_{n}^{2})},
\end{aligned}
\end{equation*}
With the two expressions above, we have $|\psi_{\tilde{T}_{n,0}}(t)-\psi_{\tilde{T}_{n,0}^{q}}(t)|{\leq}|t|[E(\tilde{T}_{n,0}-\tilde{T}_{n,0}^{q})^{2}]^{1/2}$. Note that
\begin{equation*}
\begin{aligned}
E(\tilde{T}_{n,0}-\tilde{T}_{n,0}^{q})^{2}&=E\{\sum_{r=q+1}^{kp}(\eta_{n,r}^{2}-\lambda_{n,r})/\sqrt{2tr(\Omega_{n}^{2})}\}^{2}\\
&=Var(\sum_{r=q+1}^{kp}\eta_{n,r}^{2})/\{2tr(\Omega_{n}^{2})\}\\
&~{\leq}\{\sum_{r=q+1}^{kp}\sqrt{Var(\eta_{n,r}^{2})}\}^{2}/\{2tr(\Omega_{n}^{2})\}\\
&~{\leq}(2+\Delta/n_{min})(\sum_{r=q+1}^{kp}\lambda_{n,r})^{2}/\{2tr(\Omega_{n}^{2})\}.
\end{aligned}
\end{equation*}
Therefore, we have
\begin{equation*}
\begin{aligned}
|\psi_{\tilde{T}_{n,0}}(t)-\psi_{\tilde{T}_{n,0}^{q}}(t)|&~{\leq}|t|[E(\tilde{T}_{n,0}-\tilde{T}_{n,0}^{q})^{2}]^{1/2}\\
&~{\leq}|t|(1+\Delta/2n_{min})^{1/2}(\sum_{r=q+1}^{kp}\lambda_{n,r})/\{tr(\Omega_{n}^{2})\}\\
&=|t|(1+\Delta/2n_{min})^{1/2}\sum_{r=q+1}^{kp}\varrho_{n,r}.
\end{aligned}
\end{equation*}

Let t be fixed. By Condition C4, for any fixed $q$, as $n,p\rightarrow\infty$, we have $\sum_{r=1}^{\infty}\varrho_{r}<\infty$ and
\begin{equation*}
\begin{aligned}
\sum_{r=q+1}^{kp}\varrho_{n,r}=\sum_{r=1}^{kp}\varrho_{n,r}-\sum_{r=1}^{q}\varrho_{n,r}\rightarrow\sum_{r=1}^{\infty}\varrho_{r}-\sum_{r=1}^{q}\varrho_{r}\rightarrow\sum_{r=q+1}^{\infty}\varrho_{r}.
\end{aligned}
\end{equation*}
By letting $q\rightarrow\infty$, we have $\sum_{r=q+1}^{\infty}\varrho_{r}\rightarrow0$. Therefore, for any given $\epsilon>0$, there exist $P_{1},Q_{1}$ and $N_{1}$, depending on t and $\epsilon$, for any $p{\geq}P_{1},~q{\geq}Q_{1}$ and $n{\geq}N_{1}$, we have
\begin{equation}
\begin{aligned}
|\psi_{\tilde{T}_{n,0}}(t)-\psi_{\tilde{T}_{n,0}^{q}}(t)|{\leq}\epsilon.
\end{aligned}
\end{equation}
Similarly, we have
\begin{equation*}
\begin{aligned}
&\tilde{T}_{n,0}^{q}=\frac{\sum_{r=1}^{q}(\eta_{n,r}^{2}-\lambda_{n,r})}{\sqrt{2tr(\Omega_{n}^{2})}}=\frac{\sum_{r=1}^{q}\eta_{n,r}^{2}}{\sqrt{2tr(\Omega_{n}^{2})}}-\frac{\sum_{r=1}^{q}\varrho_{n,r}}{\sqrt{2}}\\
&\tilde{T}_{n,0}^{*}=\frac{T_{n,0}^{*}-tr(\Omega_{n})}{\sqrt{2tr(\Omega_{n}^{2})}}=\frac{\sum_{r=1}^{kp}\lambda_{n,r}A_{r}-tr(\Omega_{n})}{\sqrt{2tr(\Omega_{n}^{2})}}=\frac{\sum_{r=1}^{kp}\lambda_{n,r}(A_{r}-1)}{\sqrt{2tr(\Omega_{n}^{2})}}\\
&\tilde{T}_{n,0}^{*q}=\frac{\sum_{r=1}^{q}\lambda_{n,r}(A_{r}-1)}{\sqrt{2tr(\Omega_{n}^{2})}}=\frac{\sum_{r=1}^{q}\varrho_{n,r}(A_{r}-1)}{\sqrt{2}}.
\end{aligned}
\end{equation*}
For any fixed $p{\geq}P_{1},~q{\geq}Q_{1}$, and we always have $p{\geq}q$, by the central limit theorem, when $n\rightarrow\infty$, we have $\tilde{T}_{n,0}^{q}\overset{L}{\longrightarrow}\tilde{T}_{n,0}^{*q}$ since as $n\rightarrow\infty$, $\eta_{n,r}\overset{L}{\longrightarrow}N(0,\lambda_{r})$ and $\eta_{n,r}$'s, $r=(1,\ldots,q)$ are asymptotically independent. Under Condition C3, there exists $N_{2}$, depending on p, q, t and $\epsilon$, for any $n{\geq}N_{2}$ we have
\begin{equation}
\begin{aligned}
|\psi_{\tilde{T}_{n,0}^{q}}(t)-\psi_{\tilde{T}_{n,0}^{*q}}(t)|{\leq}\epsilon.
\end{aligned}
\end{equation}
Recall that $\zeta\overset{d}{=}\sum_{r=1}^{\infty}\varrho_{r}(A_{r}-1)/\sqrt{2}$, we have $\zeta^{q}\overset{d}{=}\sum_{r=1}^{q}\varrho_{r}(A_{r}-1)/\sqrt{2}$. Under Condition C4, for any fixed q, as $p\rightarrow\infty$, we have $\tilde{T}_{n,0}^{*q}\overset{L}{\longrightarrow}\zeta^{q}$, there exists $p_{2}$, depending on q, t and $\epsilon$, for any $p{\geq}P_{2}$ we have
\begin{equation}
\begin{aligned}
|\psi_{\tilde{T}_{n,0}^{*q}}(t)-\psi_{\zeta^{q}}(t)|{\leq}\epsilon.
\end{aligned}
\end{equation}
Similarly, we have
\begin{equation*}
\begin{aligned}
|\psi_{\zeta^{q}}(t)-\psi_{\zeta}(t)|&~{\leq}|t|\left[E\left(\sum_{r=q+1}^{\infty}\varrho_{r}(A_{r}-1)/\sqrt{2}\right)^{2}\right]^{1/2}\\
&~{\leq}|t|\left[Var\left(\sum_{r=q+1}^{\infty}\varrho_{r}(A_{r}-1)/\sqrt{2}\right)\right]^{1/2}\\
&=|t|\left(\sum_{r=q+1}^{\infty}\varrho_{r}^{2}\right)^{1/2}{\leq}|t|\left(\sum_{r=q+1}^{\infty}\varrho_{r}\right),
\end{aligned}
\end{equation*}
it tends to 0 as $q{\rightarrow}\infty$ under Condition C4, there exists $Q_{2}$, depending on t and $\epsilon$, for any $q{\geq}Q_{2}$ we have
\begin{equation}
\begin{aligned}
|\psi_{\zeta^{q}}(t)-\psi_{\zeta}(t)|{\leq}\epsilon.
\end{aligned}
\end{equation}
For any $n{\geq}max(N_{1},N_{2}),~p{\geq}max(P_{1},P_{2})$ and $q{\geq}max(Q_{1},Q_{2})$, according to (A.1)-(A.4) we have
\begin{equation*}
\begin{aligned}
|\psi_{\tilde{T}_{n,0}}(t)-\psi_{\zeta}(t)|&~{\leq}|\psi_{\tilde{T}_{n,0}}(t)-\psi_{\tilde{T}_{n,0}^{q}}(t)|+|\psi_{\tilde{T}_{n,0}^{q}}(t)-\psi_{\tilde{T}_{n,0}^{*q}}(t)|\\
&+|\psi_{\tilde{T}_{n,0}^{*q}}(t)-\psi_{\zeta^{q}}(t)|+|\psi_{\zeta^{q}}(t)-\psi_{\zeta}(t)|{\leq}4\epsilon,
\end{aligned}
\end{equation*}
let $\epsilon\rightarrow0$, we have $|\psi_{\tilde{T}_{n,0}}(t)-\psi_{\zeta}(t)|\rightarrow0$, then the first expression in Theorem \ref{th2.3} (1) is proved.\\
To prove the second expression in Theorem \ref{th2.3} (1), Similarly, we have
\begin{equation*}
\begin{aligned}
|\psi_{\tilde{T}_{n,0}^{*}}(t)-\psi_{\tilde{T}_{n,0}^{*q}}(t)|&~{\leq}|t|[E(\tilde{T}_{n,0}^{*}-\tilde{T}_{n,0}^{*q})^{2}]^{1/2}\\
&=|t|\left[E\left(\sum_{r=q+1}^{p}\varrho_{p,r}(A_{r}-1)/\sqrt{2}\right)^{2}\right]^{1/2}\\
&=|t|\left(\sum_{r=q+1}^{p}\varrho_{p,r}^{2}\right)^{1/2}{\leq}|t|\sum_{r=q+1}^{p}\varrho_{p,r}.
\end{aligned}
\end{equation*}
Under Condition C4, for any given $\epsilon>0$, there exists $P_{3}$ and $Q_{3}$, depending on t and $\epsilon$, for any $p{\geq}P_{3}$ and $q{\geq}Q_{3}$ we have
\begin{equation}
\begin{aligned}
|\psi_{\tilde{T}_{n,0}^{*}}(t)-\psi_{\tilde{T}_{n,0}^{*q}}(t)|{\leq}\epsilon.
\end{aligned}
\end{equation}
For any $p{\geq}max(P_{2},P_{3})$ and $q{\geq}max(Q_{2},Q_{3})$, according to (A.3)-(A.5) we have
\begin{equation*}
\begin{aligned}
|\psi_{\tilde{T}_{n,0}^{*}}(t)-\psi_{\zeta}(t)|&~{\leq}|\psi_{\tilde{T}_{n,0}^{*}}(t)-\psi_{\tilde{T}_{n,0}^{*q}}(t)|\\
&+|\psi_{\tilde{T}_{n,0}^{*q}}(t)-\psi_{\zeta^{q}}(t)|+|\psi_{\zeta^{q}}(t)-\psi_{\zeta}(t)|{\leq}3\epsilon.
\end{aligned}
\end{equation*}
let $\epsilon\rightarrow0$, we have $|\psi_{\tilde{T}_{n,0}^{*}}(t)-\psi_{\zeta}(t)|\rightarrow0$, then the second expression in Theorem \ref{th2.3} (1) is proved.\\
Then we prove the first expression in Theorem \ref{th2.3} (2). For convenience, we denote
\begin{equation*}
\begin{aligned}
T_{n,0}=\hat{\mu}^{T}(H\otimes\W_{p})\hat{\mu}=\sum_{\alpha,\beta}c_{\alpha\beta}\bar{y}_{\alpha}^{T}W_{p}\bar{y}_{\beta}, 
\end{aligned}
\end{equation*}
where $c_{\alpha\beta}$ is the $(\alpha,\beta)$th entry of the $k{\times}k$ matrix $G^{T}G$, and $\hat{\mu}=(\bar{y}_{1}^{T},\ldots,\bar{y}_{k}^{T})^{T}$ and $\bar{y_{i}}$ is an unbiased estimate of $\mu_{i}$. Furthermore,
\begin{equation*}
\begin{aligned}
T_{n,0}&=\sum_{\alpha,\beta}c_{\alpha\beta}\bar{y}_{\alpha}^{T}W_{p}\bar{y}_{\beta}\\
&=\frac{1}{n_{\alpha}^{2}}\sum_{\alpha=1}^{k}c_{\alpha\alpha}\sum_{i,j}\bar{y}_{{\alpha}i}^{T}W_{p}\bar{y}_{{\alpha}j}+\frac{1}{n_{\alpha}n_{\beta}}\sum_{\alpha\neq\beta}c_{\alpha\beta}\sum_{i,j}\bar{y}_{{\alpha}i}^{T}W_{p}\bar{y}_{{\beta}j}\\
&=\sum_{\alpha=1}^{k}\frac{c_{\alpha\alpha}}{n_{\alpha}^{2}}\sum_{i=1}\bar{y}_{{\alpha}i}^{T}W_{p}\bar{y}_{{\alpha}i}+\sum_{\alpha=1}^{k}\frac{2c_{\alpha\alpha}}{n_{\alpha}^{2}}\sum_{i<j}\bar{y}_{{\alpha}i}^{T}W_{p}\bar{y}_{{\alpha}j}+\sum_{1\leq\alpha<\beta{\leq}k}\frac{2c_{\alpha\beta}}{n_{\alpha}n_{\beta}}\sum_{i,j}\bar{y}_{{\alpha}i}^{T}W_{p}\bar{y}_{{\beta}j},\\
\end{aligned}
\end{equation*}
thus, we have
\begin{equation*}
\begin{aligned}
E(\sum_{\alpha=1}^{k}\frac{c_{\alpha\alpha}}{n_{\alpha}^{2}}\sum_{i=1}\bar{y}_{{\alpha}i}^{T}W_{p}\bar{y}_{{\alpha}i})=\sum_{\alpha=1}^{k}\frac{c_{\alpha\alpha}}{n_{\alpha}^{2}}E(\sum_{i=1}\bar{y}_{{\alpha}i}^{T}W_{p}\bar{y}_{{\alpha}i})=\sum_{\alpha=1}^{k}\frac{c_{\alpha\alpha}}{n_{\alpha}}tr(W_{p}\Sigma_{\alpha}),
\end{aligned}
\end{equation*}
and
\begin{equation*}
\begin{aligned}
&Var(\sum_{\alpha=1}^{k}\frac{c_{\alpha\alpha}}{n_{\alpha}^{2}}\sum_{i=1}\bar{y}_{{\alpha}i}^{T}W_{p}\bar{y}_{{\alpha}i})\\
&=\sum_{\alpha=1}^{k}\frac{c_{\alpha\alpha}^{2}}{n_{\alpha}^{4}}\sum_{i=1}Var(\bar{y}_{{\alpha}i}^{T}W_{p}\bar{y}_{{\alpha}i})\\
&=\sum_{\alpha=1}^{k}\frac{c_{\alpha\alpha}^{2}}{n_{\alpha}^{4}}\sum_{i=1}\{E(\bar{y}_{{\alpha}i}^{T}W_{p}\bar{y}_{{\alpha}i}\bar{y}_{{\alpha}i}^{T}W_{p}\bar{y}_{{\alpha}i})-E^{2}(\bar{y}_{{\alpha}i}^{T}W_{p}\bar{y}_{{\alpha}i})\}\\
&=\sum_{\alpha=1}^{k}\frac{c_{\alpha\alpha}^{2}}{n_{\alpha}^{4}}\sum_{i=1}\{E(z_{{\alpha}i}^{T}\Gamma_{\alpha}^{T}W_{p}\Gamma_{\alpha}z_{{\alpha}i}z_{{\alpha}i}^{T}\Gamma_{\alpha}^{T}W_{p}\Gamma_{\alpha}z_{{\alpha}i})-E^{2}(z_{{\alpha}i}^{T}\Gamma_{\alpha}^{T}W_{p}\Gamma_{\alpha}z_{{\alpha}i})\}\\
&=\sum_{\alpha=1}^{k}\frac{c_{\alpha\alpha}^{2}}{n_{\alpha}^{4}}n_{\alpha}\{tr^{2}(\Gamma_{\alpha}^{T}W_{p}\Gamma_{\alpha})+2tr(\Gamma_{\alpha}^{T}W_{p}\Gamma_{\alpha}\Gamma_{\alpha}^{T}W_{p}\Gamma_{\alpha})+{\Delta}tr(\Gamma_{\alpha}^{T}W_{p}\Gamma_{\alpha}\circ\\
&~~~~\Gamma_{\alpha}^{T}W_{p}\Gamma_{\alpha})-tr^{2}(\Gamma_{\alpha}^{T}W_{p}\Gamma_{\alpha})\}\\
&~\leq\sum_{\alpha=1}^{k}\frac{c_{\alpha\alpha}^{2}}{n_{\alpha}^{3}}(2+\Delta)tr(W_{p}\Sigma_{\alpha})^{2},
\end{aligned}
\end{equation*}
where $\circ$ denoting the Hardmard product operator, and
\begin{equation*}
\begin{aligned}
&Var\left(\sum_{\alpha=1}^{k}\frac{2c_{\alpha\alpha}}{n_{\alpha}^{2}}\sum_{i<j}\bar{y}_{{\alpha}i}^{T}W_{p}\bar{y}_{{\alpha}j}+\sum_{1\leq\alpha<\beta{\leq}k}\frac{2c_{\alpha\beta}}{n_{\alpha}n_{\beta}}\sum_{i,j}\bar{y}_{{\alpha}i}^{T}W_{p}\bar{y}_{{\beta}j}\right)\\
&=2\left(\sum_{\alpha=1}^{k}\frac{c_{\alpha\alpha}^{2}tr(W_{p}\Sigma_{\alpha})^{2}}{n_{\alpha}^{2}}+\sum_{\alpha\neq\beta}\frac{c_{\alpha\beta}^{2}tr(W_{p}\Sigma_{\alpha}W_{p}\Sigma_{\beta})}{n_{\alpha}n_{\beta}}\right).
\end{aligned}
\end{equation*}
From the above expression, we have
\begin{equation*}
\begin{aligned}
Var(\sum_{\alpha=1}^{k}\frac{c_{\alpha\alpha}}{n_{\alpha}^{2}}\sum_{i=1}\bar{y}_{{\alpha}i}^{T}W_{p}\bar{y}_{{\alpha}i})=o\{Var(\sum_{\alpha=1}^{k}\frac{2c_{\alpha\alpha}}{n_{\alpha}^{2}}\sum_{i<j}\bar{y}_{{\alpha}i}^{T}W_{p}\bar{y}_{{\alpha}j}+\sum_{1\leq\alpha<\beta{\leq}k}\frac{2c_{\alpha\beta}}{n_{\alpha}n_{\beta}}\sum_{i,j}\bar{y}_{{\alpha}i}^{T}W_{p}\bar{y}_{{\beta}j})\}.
\end{aligned}
\end{equation*}
Furthermore, we have
\begin{equation*}
\begin{aligned}
Var(T_{n,0})&=Var(\sum_{\alpha=1}^{k}\frac{c_{\alpha\alpha}}{n_{\alpha}^{2}}\sum_{i=1}\bar{y}_{{\alpha}i}^{T}W_{p}\bar{y}_{{\alpha}i})\\
&+Var(\sum_{\alpha=1}^{k}\frac{2c_{\alpha\alpha}}{n_{\alpha}^{2}}\sum_{i<j}\bar{y}_{{\alpha}i}^{T}W_{p}\bar{y}_{{\alpha}j}+\sum_{1\leq\alpha<\beta{\leq}k}\frac{2c_{\alpha\beta}}{n_{\alpha}n_{\beta}}\sum_{i,j}\bar{y}_{{\alpha}i}^{T}W_{p}\bar{y}_{{\beta}j})\\
&+2Cov(\sum_{\alpha=1}^{k}\frac{c_{\alpha\alpha}}{n_{\alpha}^{2}}\sum_{i=1}\bar{y}_{{\alpha}i}^{T}W_{p}\bar{y}_{{\alpha}i},\sum_{\alpha=1}^{k}\frac{2c_{\alpha\alpha}}{n_{\alpha}^{2}}\sum_{i<j}\bar{y}_{{\alpha}i}^{T}W_{p}\bar{y}_{{\alpha}j}+\sum_{1\leq\alpha<\beta{\leq}k}\frac{2c_{\alpha\beta}}{n_{\alpha}n_{\beta}}\sum_{i,j}\bar{y}_{{\alpha}i}^{T}W_{p}\bar{y}_{{\beta}j}),
\end{aligned}
\end{equation*}
and denotes $\vartheta=\sum_{\alpha=1}^{k}\frac{2c_{\alpha\alpha}}{n_{\alpha}^{2}}\sum_{i<j}\bar{y}_{{\alpha}i}^{T}W_{p}\bar{y}_{{\alpha}j}+\sum_{1\leq\alpha<\beta{\leq}k}\frac{2c_{\alpha\beta}}{n_{\alpha}n_{\beta}}\sum_{i,j}\bar{y}_{{\alpha}i}^{T}W_{p}\bar{y}_{{\beta}j}$, we know
\begin{equation*}
\begin{aligned}
Cov(\sum_{\alpha=1}^{k}\frac{c_{\alpha\alpha}}{n_{\alpha}^{2}}\sum_{i=1}\bar{y}_{{\alpha}i}^{T}W_{p}\bar{y}_{{\alpha}i},\vartheta)^{2}{\leq}Var(\sum_{\alpha=1}^{k}\frac{c_{\alpha\alpha}}{n_{\alpha}^{2}}\sum_{i=1}\bar{y}_{{\alpha}i}^{T}W_{p}\bar{y}_{{\alpha}i})Var(\vartheta)=Var(\vartheta)\{1+o(1)\},
\end{aligned}
\end{equation*}
thus, we have
\begin{equation*}
\begin{aligned}
Var(T_{n,0})=Var(\vartheta)\{1+o(1)\},
\end{aligned}
\end{equation*}
and
\begin{equation*}
\begin{aligned}
\frac{\sum_{\alpha=1}^{k}\frac{c_{\alpha\alpha}}{n_{\alpha}^{2}}\sum_{i=1}\bar{y}_{{\alpha}i}^{T}W_{p}\bar{y}_{{\alpha}i}-E(\sum_{\alpha=1}^{k}\frac{c_{\alpha\alpha}}{n_{\alpha}^{2}}\sum_{i=1}\bar{y}_{{\alpha}i}^{T}W_{p}\bar{y}_{{\alpha}i})}{\sqrt{Var(\vartheta)}}=O_{p}(1).
\end{aligned}
\end{equation*}
From the first expression in Theorem \ref{th2.3} (2), we have
\begin{equation*}
\begin{aligned}
\frac{T_{n.0}-E(T_{n.0})}{\sqrt{Var(\vartheta)}}=\frac{\sum_{\alpha=1}^{k}\frac{c_{\alpha\alpha}}{n_{\alpha}^{2}}\sum_{i=1}\bar{y}_{{\alpha}i}^{T}W_{p}\bar{y}_{{\alpha}i}-E(\sum_{\alpha=1}^{k}\frac{c_{\alpha\alpha}}{n_{\alpha}^{2}}\sum_{i=1}\bar{y}_{{\alpha}i}^{T}W_{p}\bar{y}_{{\alpha}i})}{\sqrt{Var(\vartheta)}}+\frac{\vartheta-E(\vartheta)}{\sqrt{Var(\vartheta)}}.
\end{aligned}
\end{equation*}
We just need to prove $\frac{\vartheta-E(\vartheta)}{\sqrt{Var(\vartheta)}}\overset{L}{\longrightarrow}N(0,1)$. Therefore, the rest can be proved using the central limit theorem of martingale difference similar to \cite{Zhou et al:2017}.\\
Now we prove the second expression in Theorem \ref{th2.3} (2) by the Lyapunov central limit theorem. When $p\rightarrow\infty$, we have
\begin{equation*}
\begin{aligned}
tr(\Omega_{n}^{3})\sum_{r=1}^{kp}\lambda_{n,r}^{3}\leq\lambda_{n,max}\sum_{r=1}^{kp}\lambda_{n,r}^{2}=\lambda_{n,max}tr(\Omega_{n}^{2}).
\end{aligned}
\end{equation*}
and hence we have $d^{*}=\frac{tr^{3}(\Omega_{n}^{2})}{tr^{2}(\Omega_{n}^{3})}\geq[\frac{\lambda_{n,max}^{2}}{tr(\Omega_{n}^{2})}]^{-1}$, then, by the Cauchy-Schwarz inequality, we have $tr^{2}(\Omega_{n}^{3}){\leq}tr(\Omega_{n}^{4})tr(\Omega_{n}^{})$, and
\begin{equation*}
\begin{aligned}
d^{*}\geq\left[\frac{tr(\Omega_{n}^{4})}{tr^{2}(\Omega_{n}^{2})}\right]^{-1}.
\end{aligned}
\end{equation*}
Besides, we can rewrite $d^{*}=p\left[\frac{tr(\Omega_{n}^{2})}{p}\right]^{3}\left[\frac{p}{tr(\Omega_{n}^{3})}\right]^{2}$. By the condition $\lambda_{n,p,max}^{2}=o[tr(\Sigma^{2})]$, as $p\rightarrow\infty$ in \cite{Bai and Saranadasa:1996}, condition $tr(\Sigma^{4})=o[tr^{2}(\Sigma^{2})]$, as $p\rightarrow\infty$ in \cite{Chen et al:2010}, or condition $tr(\Sigma^{l})/p{\rightarrow}a_{l}\in(0,\infty),~l=(1,2,3)$, as $p\rightarrow\infty$ in \cite{Srivastava et al:2008} it follows that $p\rightarrow\infty$, $d^{*}\rightarrow\infty$. The skewness of $T_{n,0}^{*}$ is $E\{T_{n,0}^{*}-E(T_{n,0}^{*})\}^{3}/Var^{3/2}(T_{n,0}^{*})=(8/d^{*})^{1/2}\rightarrow0$, according to Lyapunov central limit theorem, we have the second expression in Theorem \ref{th2.3} (2).\\
Now, the uniform convergence result given in (\ref{eq2.11}) follows immediately from the convergence in distribution results given in Lemma 2.11 of 
\cite{van der Vaart:1998}. Denotes $\tilde{x}=[x-tr(\Omega_{n})]/[2tr(\Omega_{n}^{2})]^{1/2}$ for any real number $x$. Since the limit $\zeta$ is a continuous random variable, the expression (\ref{eq2.11}) follows directly from the convergence in distribution of both $T_{n,0}$ and $T_{n,0}^{*}$ to $\zeta$ and the triangular inequality
\begin{equation*}
\begin{aligned}
&\sup_{x}|Pr(T_{n,0}{\leq}x)-Pr(T_{n,0}^{*}{\leq}x)|\\
&=\sup_{x}|Pr(\tilde{T}_{n,0}{\leq}\tilde{x})-Pr(\tilde{T}_{n,0}^{*}{\leq}\tilde{x})|\\
&\leq\sup_{x}|Pr(\tilde{T}_{n,0}{\leq}\tilde{x})-Pr(\zeta{\leq}\tilde{x})|+\sup_{x}|Pr(\tilde{T}_{n,0}^{*}{\leq}\tilde{x})-Pr(\zeta{\leq}\tilde{x})|\\
&\rightarrow0 ~~as~~n,p\rightarrow\infty.
\end{aligned}
\end{equation*}
This completes the proof of Theorem \ref{th2.3}.\\

\textbf{Proof of Theorem \ref{th2.4}}\\
\noindent Recall the expression (\ref{eq2.14}), we have
\begin{equation*}
\begin{aligned}
\frac{R-tr(\Omega_{n})}{\sqrt{2tr(\Omega_{n}^{2})}}=\frac{\chi_{d}^{2}-d}{\sqrt{2d}}.
\end{aligned}
\end{equation*}
Next we show that when $n,p\rightarrow\infty$, $d\rightarrow\infty$ by the Cauchy–Schwarz inequality, we have
\begin{equation*}
\begin{aligned}
tr^{2}(\Omega_{n}^{3}){\leq}tr(\Omega_{n}^{2})tr(\Omega_{n}^{4}),~tr(\Omega_{n}^{3}){\geq}\frac{tr^{2}(\Omega_{n}^{2})}{tr(\Omega_{n})},
\end{aligned}
\end{equation*}
therefore, we have
\begin{equation*}
\begin{aligned}
d=\frac{tr^{2}(\Omega_{n})}{tr(\Omega_{n}^{2})}{\geq}\frac{tr^{3}(\Omega_{n}^{2})}{tr^{2}(\Omega_{n}^{3})}{\geq}\frac{tr^{2}(\Omega_{n}^{2})}{tr(\Omega_{n}^{4})}.
\end{aligned}
\end{equation*}
If we can prove that $tr(\Omega_{n}^{4})=o\{tr^{2}(\Omega_{n}^{2})\}$, we have $d\rightarrow\infty$. To prove it, we have
\begin{equation*}
\begin{aligned}
tr(\Omega_{n}^{4})=tr[\{(A{\otimes}W_{p})\Sigma\}^{4}]=\sum_{i_{1},i_{2},i_{3},i_{4}}a_{i_{1}i_{2}}a_{i_{2}i_{3}}a_{i_{3}i_{4}}a_{i_{4}i_{1}}tr(W_{p}\Sigma_{i_{1}}W_{p}\Sigma_{i_{2}}W_{p}\Sigma_{i_{3}}W_{p}\Sigma_{i_{4}}).
\end{aligned}
\end{equation*}
Then, by Condition C5 we have for $i_{1},i_{2},i_{3},i_{4}\in\{1,\ldots,k\}$, as $p\rightarrow\infty$
\begin{equation*}
\begin{aligned}
tr(W_{p}\Sigma_{i1}W_{p}\Sigma_{i2}W_{p}\Sigma_{i3}W_{p}\Sigma_{i4})=o\{tr(W_{p}\Sigma_{i1}W_{p}\Sigma_{i2})tr(W_{p}\Sigma_{i3}W_{p}\Sigma_{i4})\}.
\end{aligned}
\end{equation*}
Then by the Cauchy–Schwarz inequality, we have $a_{ij}=h_{ij}/\sqrt{n_{i}}\sqrt{n_{j}}=g_{i}^{T}g_{j}/\sqrt{n_{i}}\sqrt{n_{j}}\leq\sqrt{a_{ii}a_{jj}},~a_{i_{1}i_{2}}a_{i_{2}i_{3}}a_{i_{3}i_{4}}a_{i_{4}i_{1}}{\leq}a_{i_{1}i_{1}}a_{i_{2}i_{2}}a_{i_{3}i_{3}}a_{i_{4}i_{4}}$, where $H:(h_{ij})_{i,j}^{k}=G^{T}G$, and $G=(g_{1},\ldots,g_{k})$. Thus, we have
\begin{equation*}
\begin{aligned}
tr(\Omega_{n}^{4})=o\left\{\sum_{i_{1}i_{2}}a_{i_{1}i_{1}}a_{i_{2}i_{2}}tr(W_{p}\Sigma_{i_{1}}W_{p}\Sigma_{i_{2}})\sum_{i_{3}i_{4}}a_{i_{3}i_{3}}a_{i_{4}i_{4}}tr(W_{p}\Sigma_{i_{3}}W_{p}\Sigma_{i_{4}})\right\}=o\{tr^{2}(\Omega_{n}^{2})\}.
\end{aligned}
\end{equation*}
Therefore, we have $d\rightarrow\infty$ and the expression (\ref{eq2.14}) is established. The second expression in Theorem \ref{th2.4} can be proved by the same methods similar to that in expression (\ref{eq2.11}) of Theorem \ref{th2.3}. This completes the proof of Theorem \ref{th2.4}.\\

\textbf{Proof of Theorem \ref{th2.5}}\\
\noindent According to the expression (\ref{eq2.16}), if we want to get the ratio-consistent estimators of $tr(\Omega_{n}),~tr^{2}(\Omega_{n})$ and $tr(\Omega_{n}^{2})$, it is equivalent to get the ratio-consistent estimators of $tr(W_{p}\Sigma_{i}),~tr^{2}(W_{p}\Sigma_{i}),~tr\{(W_{p}\Sigma_{i})^{2}\}$ and $tr(W_{p}\Sigma_{i})tr(W_{p}\Sigma_{j}),~tr(W_{p}\Sigma_{i}W_{p}\Sigma_{j})$. Under Condition C1-C3, as $n_{i}\rightarrow\infty$, similar to Lemma S.3 in the supplementary material in \cite{Zhang:2020}, the ratio-consistent estimators of $tr(W_{p}\Sigma_{i}),~tr^{2}(W_{p}\Sigma_{i}),~tr\{(W_{p}\Sigma_{i})^{2}\}$ is shown in the expression (\ref{eq2.17}), it follows that under Condition C1-C3, as $n_{i}\rightarrow\infty$, $tr(W_{p}\hat{\Sigma}_{i})tr(W_{p}\hat{\Sigma}_{j})$ is also ratio-consistent for $tr(W_{p}\Sigma_{i})tr(W_{p}\Sigma_{j})$ uniformly for all p, where $\Sigma_{i}$ is defined in the preceding text. Besides, similar to the proof of Theorem 2 of \cite{Zhang:2021}, we get the ratio-consistent estimators of $tr(W_{p}\Sigma_{i}W_{p}\Sigma_{j})$ is $tr(W_{p}\hat{\Sigma}_{i}W_{p}\hat{\Sigma}_{j})$.\\
Therefore, under Condition C1-C3, as $n_{i}\rightarrow\infty$, we have
\begin{equation*}
\begin{aligned}
tr(\widehat{\Omega}_{n})=\sum_{i=1}^{k}a_{ii}tr(W_{p}\Sigma_{i})\{tr(W_{p}\hat{\Sigma}_{i})/tr(W_{p}\Sigma_{i})\}=tr(\Omega_{n})\{1+o_{p}(1)\}
\end{aligned}
\end{equation*}
uniformly for all p, thus, we have $tr(\widehat{\Omega}_{n})/tr(\Omega_{n})\overset{P}{\longrightarrow}1$ uniformly for all p. similarly, we can show that $\widehat{tr^{2}(\Omega_{n})}/tr^{2}(\Omega_{n})\overset{P}{\longrightarrow}1$ and $\widehat{tr(\Omega_{n}^{2})}/tr(\Omega_{n}^{2})\overset{P}{\longrightarrow}1$ under the same conditions by the expression (\ref{eq2.16}) and (\ref{eq2.18}). It follows that under Condition C1-C3, as $n_{i}\rightarrow\infty$, we have
\begin{equation*}
\begin{aligned}
\frac{\hat{\beta}}{\beta}=\frac{\widehat{tr(\Omega_{n}^{2})}/tr(\Omega_{n}^{2})}{tr(\widehat{\Omega}_{n})/tr(\Omega_{n})}\overset{P}{\longrightarrow}1,~\frac{\hat{d}}{d}=\frac{\widehat{tr^{2}(\Omega_{n})}/tr^{2}(\Omega_{n})}{\widehat{tr(\Omega_{n}^{2})}/tr(\Omega_{n}^{2})}\overset{P}{\longrightarrow}1,
\end{aligned}
\end{equation*}
uniformly for all p. Therefore, we have
\begin{equation*}
\begin{aligned}
\frac{\widehat{\beta}\chi_{\widehat{d}}^{2}(\alpha)}{\beta\chi_{d}^{2}(\alpha)}\overset{P}{\longrightarrow}1.
\end{aligned}
\end{equation*}
uniformly for all p. This completes the proof of Theorem \ref{th2.5}.\\

\textbf{Proof of Theorem \ref{th2.6}}\\
\noindent Under the expression (\ref{eq2.23}) and (\ref{eq2.4}), we have
\begin{equation*}
\begin{aligned}
T_{n}=\{T_{n,0}+\mu^{T}(H\otimes\W_{p})\mu\}\{1+o_{p}(1)\}=\{T_{n,0}+tr(W_{p}M^{T}HM)\{1+o_{p}(1)\},
\end{aligned}
\end{equation*}
where denotes that $\mu^{T}(H\otimes\W_{p})\mu\}=tr(W_{p}M^{T}HM)$. Thus, we have
\begin{equation*}
\begin{aligned}
&Pr\left\{T_{n}>\widehat{\beta}\chi_{\widehat{d}}^{2}(\alpha)\right\}\\
&=Pr\left\{T_{n,0}-tr(\Omega_{n})\geq\widehat{\beta}\chi_{\widehat{d}}^{2}(\alpha)-tr(\Omega_{n})-tr(W_{p}M^{T}HM)\right\}\{1+o(1)\}\\
&=Pr\left\{\frac{T_{n,0}-tr(\Omega_{n})}{\sqrt{2tr(\Omega_{n}^{2})}}\geq\frac{\widehat{\beta}\chi_{\widehat{d}}^{2}(\alpha)-tr(\Omega_{n})}{\sqrt{2tr(\Omega_{n}^{2})}}-\frac{ntr\{W_{p}M^{T}(n^{-1}H)M\}}{\sqrt{2tr(\Omega_{n}^{2})}}\right\}\{1+o(1)\}.
\end{aligned}
\end{equation*}
Note that $H^{*}=\lim_{n\rightarrow\infty}n^{-1}H$, under Condition C1-C4, Theorem \ref{th2.3} (1) and Theorem \ref{th2.5} lead to
\begin{equation*}
\begin{aligned}
&Pr\left\{T_{n}>\widehat{\beta}\chi_{\widehat{d}}^{2}(\alpha)\right\}\\
&=Pr\left\{\zeta\geq\frac{\widehat{\beta}\chi_{\widehat{d}}^{2}(\alpha)-tr(\Omega_{n})}{\sqrt{2tr(\Omega_{n}^{2})}}-\frac{ntr\{W_{p}M^{T}H^{*}M\}}{\sqrt{2tr(\Omega_{n}^{2})}}\right\}\{1+o(1)\}\\
&=Pr\left\{\zeta\geq\frac{\chi_{\widehat{d}}^{2}(\alpha)-d}{\sqrt{2d}}-\frac{ntr\{W_{p}M^{T}H^{*}M\}}{\sqrt{2tr(\Omega_{n}^{2})}}\right\}\{1+o(1)\}
\end{aligned}
\end{equation*}
Similarly, under Condition C1-C3 and C5, Theorem \ref{th2.3} (2) and Theorem \ref{th2.5} lead to
\begin{equation*}
\begin{aligned}
Pr\left\{T_{n}>\widehat{\beta}\chi_{\widehat{d}}^{2}(\alpha)\right\}=\Phi\left\{-z_{\alpha}+\frac{ntr\{W_{p}M^{T}H^{*}M\}}{\sqrt{2tr(\Omega_{n}^{2})}}\right\}\{1+o(1)\}.
\end{aligned}
\end{equation*}
This completes the proof of Theorem \ref{th2.6}.\\

\end{appendices}

\end{document}